\documentclass[11pt,twoside,leqno,mathscr,euscript]{aomamlt2e} 
  \pageno{993}
 \renewcommand{\theequation}{\thesection.\arabic{equation}}
\received{June 3, 2004}
 \def\eref#1{{\rm (\ref{#1})}}

\def\mop#1{\mathop{\symb{\rm#1}}\nolimits}
\def\moplim#1{\mathop{\hbox{\rm#1}}}
\let\emptyset=\phi
\let\phi=\varphi
\let\orho\rho 
\def\rho{\varrho}
\def\eps{\varepsilon}

\def\d{\partial}

\newcommand\hfq{\hfill\qed}
  \newcommand{\ee}{{\hskip1pt\rm \'{\hskip-6.5pt \it e}}}

\def\({\left(}
\def\){\right)}

\def\newop#1{\expandafter\def\csname#1\endcsname{\mop{#1}}}
\def\newoplim#1{\expandafter\def\csname#1\endcsname{\moplim{#1}}}
\newop{deg}
\newop{div}
\newop{supp}
\newop{sin}
\newop{arcsin}
\newop{rank}
\newop{tr}
\newop{cos}
\newop{arccos}
\newop{tan}
\newop{sinh}
\newop{arsinh}
\newop{cosh}
\newop{arcosh}
\newop{tanh}
\newop{arctan}
\newop{artanh}
\newop{log}
\newop{exp}
\newop{sgn}
\newop{dim}
\newop{ker}
\newoplim{min}
\newoplim{max}
\newoplim{sup}
\newoplim{inf}
\newoplim{lim}

\def\limsup{\moplim{lim sup}}

\def\CE{{\cal E}}

\def\CU{{\cal U}}

\def\CP{{\cal P}}

\def\CF{{\cal F}}
\def\CG{{\cal G}}
\def\CH{{\cal H}}

\def\CZ{{\cal Z}}
\def\CX{{\cal X}}
\def\CC{{\cal C}}

\def\CB{{\cal B}}

\newcommand\cinf[1][ ]{\ifthenelse{\equal{#1}{ }}{{\cal C}^\infty}{{\cal C}^\infty(#1)}}
\newcommand\cOinf[1][ ]{\ifthenelse{\equal{#1}{ }}{{\cal C}_0^\infty}{{\cal C}_0^\infty(#1)}}
\newcommand\symb[2][\bf]{{\mathchoice{\hbox{#1#2}}{\hbox{#1#2}}%
        {\hbox{\scriptsize#1#2}}{\hbox{\tiny#1#2}}}}

\def\L^#1{{\rm L}^{\!#1}}
\def\R{{\symb R}}
\def\N{{\symb N}}
\def\Z{{\symb Z}}

\newcommand\scal[1]{\left\langle{#1}\right\rangle}
\newtheorem{lemma}{Lemma}[section]
\newtheorem{theorem}[lemma]{Theorem}
\newtheorem{proposition}[lemma]{Proposition}
\newtheorem{corollary}[lemma]{Corollary}

\newcommand\prop[1]{Proposition~{\rm \ref{#1}}}
\newcommand\lem[1]{Lemma~{\rm \ref{#1}}}
\newcommand\theo[1]{Theorem~{\rm \ref{#1}}}
\newcommand\cor[1]{Corollary~{\rm \ref{#1}}}

\theorembodyfont{\upshape}
\newtheorem{example}[lemma]{{\it Example}}

\newtheorem{definition}[lemma]{{\it Definition}}
\newtheorem{remark}[lemma]{{\it Remark}}

 \def\expect{{\mathbf E}}
\def\TV{{\mathrm{TV}}}  
\def\cC{{\mathscr{C}}} \def\cD{{\mathscr{D}}} 
\def\qsf/{asymptotically strong
Feller} \def\qsfp/{asymptotic strong Feller property} \def\X{{\cal X}}
 \def\P{{\mathbf P}} \def\J{{J}} \def\K{{K}}
\def\KK{{\cal K}} 
\def\Ltwo{\L^2_{0}}

\newcommand{\eqdef}{\stackrel{\mbox{\tiny\rm def}}{=}}

\newcommand{\TT}{\mathbb{T}}
\newcommand{\domain}{\TT^2}
\newcommand{\EE}{\expect}
\newcommand{\tM}{\widetilde{M}}

\newcommand{\norm}[1]{\vert\!\vert\!\vert #1 \vert\!\vert\!\vert}
\newcommand{\uu}{\xi}
\newcommand{\err}{\orho}
\newcommand{\Jb}{\check\J}
\newcommand{\Ja}{\hat \J}
\newcommand{\A}{A}
\newcommand{\M}{M}
\newcommand{\Ball}{\mathcal{B}}

\begin{document}
\currannalsline{164}{2006} 

 \title{Ergodicity of the\\ 2D Navier-Stokes equations\\ with degenerate
stochastic forcing}

\twoauthors{Martin Hairer}{Jonathan C.\ Mattingly}

 \institution{The University of Warwick, Coventry, United Kingdom\\
\email{hairer@maths.warwick.co.uk} \\ \vglue-9pt
 Duke University, Durham NC\\
\email{jonm@math.duke.edu}
\\
\phantom{down}}

 \shorttitle{Ergodicity of the 2D Navier-Stokes equations}  

\centerline{\bf Abstract}\vglue12pt
  The stochastic 2D Navier-Stokes equations on the torus driven by
  degenerate noise are studied. We characterize the smallest closed
  invariant subspace for this model and show that the dynamics
  restricted to that subspace is ergodic.  In particular, our results
  yield a purely geometric characterization of a class of noises for
  which the equation is ergodic in $\L^2_0(\TT^2)$. Unlike previous
  works, this class is independent of the viscosity and the strength
  of the noise. The two main tools of our analysis are the
  \textit{asymptotic strong Feller} property, introduced in this work,
  and an approximate integration by parts formula. The first, when
  combined with a weak type of irreducibility, is shown to 
  ensure that the dynamics is ergodic. The second is used to show
  that the first holds under a H{\"o}rmander-type condition. This
  requires some interesting nonadapted stochastic analysis.

\section{Introduction}

In this article, we investigate the ergodic properties of the 2D Navier-Stokes equations.
Recall that the Navier-Stokes equations describe the time evolution of an incompressible fluid
and are given by
\begin{equation} \label{e:SNS}
\d_t u  + (u \cdot \nabla) u = \nu \Delta u - \nabla p + \xi\;,\quad
{\rm div} \ u = 0\;,
\end{equation}   
where $u(x,t) \in \R^2$ denotes the value of the velocity field at time $t$ and position $x$, $p(x,t)$
denotes the pressure, and $\xi(x,t)$ is an external force field acting on the fluid.
We will consider the case when $x \in \TT^2$, the two-dimensional torus. 
Our mathematical model for the driving force $\xi$ is a Gaussian field which is
white in time and colored in space. We are particularly interested in the case when only a
few Fourier modes of $\xi$ are nonzero, so that there is a well-defined ``injection scale''
$L$ at which energy is pumped into the system. 
Remember that both the energy $\|u\|^2 = \int |u(x)|^2\,dx$ and the enstrophy
$\|\nabla \wedge u\|^2$ are invariant under the nonlinearity of the 2D 
Navier-Stokes equations (i.e.\  they are preserved by the flow of \eref{e:SNS} if $\nu = 0$ and $\xi = 0$).

From a careful study of the nonlinearity (see e.g.\   \cite{MR2003k:76071} for a survey and \cite{MR1914189} for some mathematical results in this field),
one expects the enstrophy to cascade down to smaller and smaller scales, until it reaches a
``dissipative scale'' $\eta$ at which the viscous term $\nu \Delta u$ dominates the 
nonlinearity $(u \cdot \nabla) u$ in \eref{e:SNS}. This picture is
complemented by that of an inverse cascade of the energy towards
larger and larger scales, until it is dissipated by finite-size effects as it reaches scales of order one.
The physically interesting range of parameters for \eref{e:SNS}, where one expects to see
both cascades and where the behavior
of the solutions is dominated by the nonlinearity, thus corresponds~to
\begin{equation} \label{e:scales}
1 \ll L^{-1} \ll \eta^{-1}\;.
\end{equation}   
The main assumptions usually made in the physics literature when discussing the
behavior of \eref{e:SNS} in the turbulent regime are ergodicity and statistical translational 
invariance of the stationary state.
We give a simple geometric characterization of a class of forcings for which \eref{e:SNS}
is ergodic, including a forcing that acts only on $4$ degrees of freedom ($2$ Fourier modes). This characterization is independent 
of the viscosity and  is shown to be sharp in a certain sense. In particular, it covers the range of
parameters \eref{e:scales}. Since we show that the invariant
measure for \eref{e:SNS} is unique, its translational invariance follows immediately from the
translational invariance of the equations.

From the mathematical point of view, the ergodic properties for infinite-dimen\-sional
systems are a field that has been intensely studied over the past two decades but is yet
 in its infancy  compared to the corresponding theory for finite-dimensional
systems. In particular, there is a gaping lack of results for truly hypoelliptic nonlinear systems, where
the noise is transmitted to the relevant degrees of freedom only through the drift.
The present article is an attempt to close this gap, at least for the particular case
of the 2D Navier-Stokes equations.
This particular case (and some closely related problems) has been an intense subject of study in recent
years. However the results obtained so far require either a
nondegenerate forcing on the ``unstable'' part of the equation
\cite{b:EMattinglySinai00},
\cite{b:KuksinShirikyan00},
\cite{b:BricmontKupiainenLefevere01},
\cite{ArmenKuk},
\cite{MatNS},
\cite{b:BricmontKupiainenLefevere02},
\cite{HExp02},
\cite{b:MasmoudiYoung02},
or the strong Feller property to hold. The latter was obtained only when the forcing acts on an infinite number
of modes
\cite{FM},
\cite{b:Fe97},
\cite{EH3},
\cite{b:MattinglySuidan05}. The former used a change
of measure via Girsanov's theorem and  the
pathwise contractive properties of the dynamics to prove ergodicity.
In all of these works, 
the noise was sufficiently nondegenerate to allow in a way for an adapted analysis (see Section~\ref{sec:EssEllip} below for the meaning of ``adapted''
in this context). 

We give a fairly complete analysis of the conditions needed to ensure
the ergodicity of the two dimensional Navier-Stokes equations. To do
so, we employ information on the structure of the nonlinearity from
\cite{b:EMattingly00} which was developed there to prove ergodicity of
the finite dimensional Galerkin approximations under conditions on the
forcing similar to this paper. However, our approach to the full PDE
is necessarily different and informed by the pathwise contractive
properties and high/low mode splitting explained in the stochastic
setting in \cite{b:Mattingly98b},
\cite{b:Mattingly98} and the ideas of
determining modes, inertial manifolds, and invariant subspaces in
general from the deterministic PDE literature (cf.
\cite{b:FoiasProdi67},
\cite{b:CoFo88}). More directly, this paper builds on
the use of the high/low splitting to prove ergodicity as first
accomplished contemporaneously in
\cite{b:BricmontKupiainenLefevere01},
\cite{b:EMattinglySinai00},
\cite{b:KuksinShirikyan00}
in the ``essentially elliptic'' setting (see section
\ref{sec:EssEllip}). In particular, this paper is the culmination of a
sequence of papers by the authors and their collaborators
\cite{b:Mattingly98b},
\cite{b:Mattingly98},
\cite{EH3},
\cite{b:EMattinglySinai00},
\cite{MatNS,HExp02},
\cite{b:Mattingly03Pre}
using these and related ideas to prove ergodicity. Yet, this is the
first to prove ergodicity of a stochastic PDE in a hypoelliptic
setting under conditions which compare favorably to those under which
similar theorems are proven for finite dimensional stochastic
differential equations. One of the keys to accomplishing this is a
recent result from \cite{b:MattinglyPardoux03Pre} on the regularity of
the Malliavin matrix in this setting.

One of the main technical contributions of the present work is to
provide an infinitesimal replacement for Girsanov's theorem in the
infinite dimensional nonadapted setting which the application of
these ideas to the fully hypoelliptic setting seems to require.
Another of the principal technical contributions is to observe that
the strong Feller property is neither essential nor natural for the
study of ergodicity in dissipative infinite-dimensional systems and to
provide an alternative. We define instead a weaker \textit{asymptotic
  strong Feller} property which is satisfied by the system under
consideration and is sufficient to give ergodicity. In many
dissipative systems, including the stochastic Navier-Stokes equations,
only a finite number of modes are unstable. Conceivably, these systems
are ergodic even if the noise is transmitted only to those unstable
modes rather than to the whole system.  The asymptotic strong Feller
property captures this idea. It is sensitive to the
regularization of the transition densities due to both probabilistic
and dynamic mechanisms.

This paper is organized as follows. In Section~\ref{sec:basicSettingMainResults}
 the precise mathematical formulation of the problem and the main results
for the stochastic Navier-Stokes equations are given. In Section~\ref{sec:abstract} we define the
\qsfp/ and prove in \theo{theo:quasiStrongFeller} that, together with an irreducibility property it
implies ergodicity of the system. 
We thus obtain the analog in our setting of the classical result often derived from
theorems of Khasminskii and Doob which states that topological irreducibility,
together with the strong Feller property, implies uniqueness of the invariant measure. 
The main technical results are given in Section~\ref{sec:SNS}, where we show how to
apply the abstract results to our problem. Although this section is written with the
stochastic Navier-Stokes equations in mind,  most of the corresponding results 
hold for a much wider class of stochastic PDEs with polynomial nonlinearities.

\demo{Acknowledgements}
  We would like to thank G. Ben Arous, W. E. J.\ Hanke, X.-M.\ Li, E.\ Pardoux,
  M.\ Romito and Y.\ Sinai for motivating and useful discussions. We would also
  like to thank the anonymous referees for their careful reading of
  the text and their subsequent corrections and useful suggestions.
  The work of MH is partially supported by the Fonds National Suisse.
  The work of JCM was partially supported by the Institut
  Universitaire de France.

\section{Setup and main results}
\label{sec:basicSettingMainResults}

Consider the two-dimensional, incompressible
Navier-Stokes equations on the to\-rus $\domain = [-\pi,\pi]^2$ driven by
a degenerate noise. Since the velocity and vorticity formulations
are equivalent in this setting, we choose to use the vorticity equation as this
simplifies the exposition. For $u$ a divergence-free velocity
field, we define the vorticity $w$ by $w = \nabla \wedge u = \d_2 u_1 - \d_1 u_2$.
Note that $u$ can be recovered from $w$ and the condition $\nabla \cdot u = 0$.
With this notation the vorticity formulation for the stochastic Navier-Stokes equations
is as follows:  
\begin{equation} \label{e:vort}
dw = \nu \Delta w\,dt + B(\KK w,w)\,dt + Q\,dW(t)\;,
\end{equation}   
where $\Delta$ is the Laplacian with periodic boundary conditions and
$B(u,w) = -(u\cdot\nabla)w$, the usual Navier-Stokes nonlinearity. The symbol
$Q\,dW(t)$ denotes a Gaussian noise process which is white in time and whose spatial
correlation structure will be described later.
The operator $\KK$ is defined in Fourier space by $(\KK w)_k = -i w_k
k^{\perp}/\|k\|^2$, where $(k_1,k_2)^\perp = (k_2,-k_1)$. By
$w_k$, we mean the scalar product of $w$ with $(2\pi)^{-1}\exp(i k\cdot x)$. It has the
property that the divergence of $\KK w$ vanishes and that $w = \nabla
\wedge (\KK w)$. Unless otherwise stated, we consider \eref{e:vort} as an equation
in $\CH = \Ltwo$, the space of real-valued
square-inte\-gra\-ble functions on the torus with vanishing mean.
Before we go on to describe the noise process $Q W$, it is
instructive to write down the two-dimensional Navier-Stokes equations (without
noise) in Fourier space:
\begin{equation} \label{e:vortFourier}
\dot w_k = - \nu |k|^2 w_k - {1\over 4\pi} \! \sum_{j + \ell = k} \scal{j^\perp\!\!\!,\ell} \Bigl({1\over |\ell|^{2}} - {1\over |j|^{2}}\Bigr)
w_j w_\ell\;.
\end{equation}   
From \eref{e:vortFourier}, we see clearly that any closed subspace of $\CH$ spanned by Fourier modes
corresponding to a subgroup of $\Z^2$ is invariant under the dynamics. In other
words, if the initial condition has a certain type of periodicity, it will be retained
by the solution for all times.

In order to describe the noise $Q\,dW(t)$, we start by introducing a convenient
way to index the Fourier basis of $\CH$. We write $\Z^2 \setminus \{(0,0)\} = \Z^2_+ \cup \Z^2_-$,
where
\begin{eqnarray*}
\Z^2_+ &= &\bigl\{(k_1,k_2) \in \Z^2\,|\, k_2 > 0\bigr\} \cup \bigl\{(k_1,0) \in \Z^2\,|\, k_1 > 0\bigr\}\;,\\
\Z^2_- &=& \bigl\{(k_1,k_2) \in \Z^2\,|\, -k \in \Z^2_+\bigr\}\;,
\end{eqnarray*}
(note that  $\Z^2_+$ is essentially the upper half-plane)
and set, for $k \in \Z^2 \setminus \{(0,0)\}$,
\begin{equation}
f_k(x) = \left\{\begin{array}{rl} \sin(k\cdot x) & \text{if $k \in \Z^2_+$,} \\[0.5em] 
\cos(k\cdot x) & \text{if $k\in\Z^2_-$.} \end{array}\right.
\end{equation}   
We also fix a set
\begin{equation} \label{e:defZ0}
\CZ_0 = \{k_n\,|\, n=1,\ldots,m\} \subset \Z^2 \setminus \{(0,0)\}\;,
\end{equation}   
which encodes the geometry of the driving noise. The set $\CZ_0$ will correspond to
the set of driven modes of equation \eref{e:vort}.

The process $W(t)$ is an $m$-dimensional Wiener process on a probability space
$(\Omega,\CF,\P)$. For definiteness, we choose $\Omega$ to be the Wiener space 
$\CC_0([0,\infty),\R^m)$,  $W$ the canonical
process, and $\P$ the Wiener measure. We denote expectations
with respect to $\P$ by $\EE$ and define $\CF_t$ to be the
$\sigma$-algebra generated by the increments of $W$ up to time $t$. 
We also denote by $\{e_n\}$ the canonical basis of $\R^m$.
The linear map $Q:\R^m \to \CH$ is given by $Q
e_n = q_n f_{k_n}$, where the $q_n$ are some strictly positive numbers, and the wave numbers
$k_n$ are given by the elements of $\CZ_0$. 
With these definitions, $QW$ is an $\CH$-valued Wiener process. We also denote the 
average rate at which energy is injected into our system by
$\mathcal{E}_0 = \tr QQ^*=\sum_n q_n^2$.

We assume that
the set $\CZ_0$ is symmetric, i.e.\  that if $k \in \CZ_0$, then 
$-k \in \CZ_0$. This is not a strong restriction and is made only to simplify the
statements of our results. It also helps to avoid the possible confusion arising from the
slightly nonstandard definition of the basis $f_k$. 
This assumption always holds for example if the noise process 
$QW$ is taken to be translation invariant.
In fact, \theo{theo:realmain} below
holds for nonsymmetric sets $\CZ_0$ if one replaces $\CZ_0$ in the
 theorem's conditions by its symmetric part.

It is well-known \cite{b:Fl94}, \cite{b:MikuleviciusRozovskii02} that
\eref{e:vort} defines a stochastic flow on $\CH$. By a stochastic flow, 
we mean a family of continuous maps $\Phi_t\colon\Omega \times
\CH \rightarrow \CH$ such that $w_t=\Phi_t(W, w_0)$ 
is the solution to \eref{e:vort} with initial condition $w_0$ and noise $W$. 
Hence, its transition semigroup $\CP_t$ given by
$\CP_t\phi(w_0) = \EE_{w_0}\phi(w_t)$ is Feller. Here, $\phi$ denotes any bounded measurable function from $\CH$ to $\R$ and we use the notation $\EE_{w_0}$ for expectations with
respect to solutions to \eref{e:vort} with initial condition $w_0$. Recall that
an \textit{invariant measure} for \eref{e:vort} is a probability measure $\mu_\star$ on $\CH$ such
that $\CP_t^*\mu_\star = \mu_\star$, where $\CP_t^*$ is the semigroup on measures dual to $\CP_t$.
While the existence of an invariant measure for \eref{e:vort} can be proved by
``soft'' techniques using the regularizing and dissipativity properties of the flow 
\cite{MR90c:35161}, \cite{b:Fl94}, showing its uniqueness is a challenging problem that requires a 
detailed analysis of the nonlinearity. The importance of showing the uniqueness of $\mu_\star$
is illustrated by the fact that it implies
\begin{equation}
\lim_{T \to \infty} {1\over T} \int_0^T \phi(w_t)\,dt = \int_{\CH} \phi(w)\,\mu_\star(dw)\;,
\end{equation}   
for all bounded continuous functions $\phi$ and $\mu_\star$-almost every initial condition $w_0 \in \CH$. It thus
gives some mathematical ground to the \textit{ergodic assumption} usually made in the
physics literature in a  discusion of the qualitative behavior of \eref{e:vort}.
The main results of this article are summarized by the following theorem:

\begin{theorem}\label{theo:realmain}
Let $\CZ_0$ satisfy the following two assumptions\/{\rm :}\/
\begin{itemize}
\ritem{\bf A1.} There exist at least two elements in $\CZ_0$ with different Euclidean norms.
\ritem{\bf A2.}  Integer linear combinations of elements of $\CZ_0$ generate $\Z^2$.
\end{itemize}
Then{\rm ,} \eref{e:vort} has a unique invariant measure in $\CH$.
\end{theorem} 

\begin{remark}
  As pointed out by J. Hanke, condition \textbf{A2} above is equivalent to
  the easily verifiable condition that the greatest common divisor of
  the set $\big\{ \det(k,\ell)
  : k,\ell \in \CZ_0 \big\}$ is $1$, where $\det(k,\ell)$ is the
  determinant of the $2 \times 2$ matrix with columns $k$ and $\ell$. 
\end{remark}

The proof of \theo{theo:realmain} is given by combining Corollary~\ref{c:ergodic} with
Proposition~\ref{prop:gen} below.
A partial converse of this ergodicity result is given by the following theorem, which
is an immediate consequence of Proposition~\ref{prop:gen}.

\begin{theorem}\label{theo:converse}
There are two qualitatively different ways in which the hypotheses of \theo{theo:realmain}
can fail. In each case there is a unique invariant measure supported on $\tilde \CH${\rm ,}
the smallest closed linear subspace of $\CH$ which is invariant under \eref{e:vort}.
\begin{itemize}
\item In the first case the elements of $\CZ_0$ are all collinear or of the same Euclidean length.
 Then $\tilde\CH$ is the finite-dimensional space spanned by $\{f_k\,|\, k \in \CZ_0\}${\rm ,} 
 and the dynamics restricted to $\tilde\CH$ is that of an Ornstein-Uhlenbeck process.
\item In the second case let $\CG$ be the smallest subgroup of $\Z^2$ containing $\CZ_0$. Then
$\tilde\CH$ is the space spanned by $\{f_k\,|\, k \in \CG\setminus\{(0,0)\}\}$. 
Let $k_1${\rm ,} $k_2$ be two generators for $\CG$ and define $v_i = 2\pi k_i / |k_i|^2${\rm ,} then $\tilde \CH$ is 
the space of functions that are periodic with respect to the translations $v_1$ and $v_2$.
\end{itemize}
\end{theorem}

\begin{remark}
That $\tilde \CH$ constructed above is invariant is clear; that it is the smallest invariant subspace
follows from the fact  that the transition probabilities of \eref{e:vort} have
a density with respect to the Lebesgue measure when projected onto any finite-dimensional
subspace of $\tilde\CH$; see \cite{b:MattinglyPardoux03Pre}.

By \theo{theo:converse} if the conditions of \theo{theo:realmain} 
are not satisfied then one of the modes with lowest wavenumber is in $\tilde\CH^\perp$. 
In fact either $f_{(1,0)} \perp \tilde\CH$ or $f_{(1,1)} \perp \tilde\CH$. On the other hand 
for sufficiently small values of $\nu$ the low modes of
\eref{e:vort} are expected to be linearly unstable \cite{b:Frisch95}.
If this is the case, a solution to \eref{e:vort} starting in $\tilde\CH^\perp$ will not converge to $\tilde\CH$ and \eref{e:vort} is therefore expected to have several distinct invariant measures on $\CH$.
It is however known that the invariant measure is unique if the viscosity is sufficiently high;
see \cite{b:Mattingly98}. (At high viscosity, all modes are linearly
stable. See \cite{b:Mattingly03Pre} for a more streamlined presentation.)
\end{remark}
  
\begin{example}
The set $\CZ_0 = \{(1,0),(-1,0),(1,1),(-1,-1)\}$ satisfies the assumptions of \theo{theo:realmain}. 
Therefore, \eref{e:vort} with noise given by 
\begin{eqnarray*}
   QW(t,x)&= &W_1(t)\sin x_1+ W_2(t)\cos x_1 +
   W_3(t)\sin(x_1+x_2) \\ &&  + W_4(t)\cos(x_1+x_2) \;,
\end{eqnarray*}
has a unique invariant measure in $\CH$ for every value of the viscosity $\nu > 0$.
\end{example}

\begin{example}
Take $\CZ_0 = \{(1,0),(-1,0),(0,1),(0,-1)\}$ whose elements are of length $1$.
Therefore, \eref{e:vort} with noise given by 
\begin{equation}
   QW(t,x) = W_1(t)\sin x_1 + W_2(t)\cos x_1  +
   W_3(t)\sin x_2 + W_4(t)\cos x_2  \;,
\end{equation}   
reduces to an Ornstein-Uhlenbeck process on the space spanned by $\sin x_1$, $\cos x_1$, $\sin x_2$, and $\cos x_2$.
\end{example}

\begin{example}
Take $\CZ_0 = \{(2,0),(-2,0),(2,2),(-2,-2)\}$, which corresponds to case 2 of \theo{theo:converse}
with $\CG$ generated by $(0,2)$ and $(2,0)$. In this case, $\tilde \CH$ is the set of functions that 
are $\pi$-periodic in both arguments.
Via the change of variables $x \mapsto x/2$, one can easily see from \theo{theo:realmain}
that \eref{e:vort} then has a unique invariant measure on $\tilde\CH$ (but not
necessarily on $\CH$).
\end{example}


\section{An abstract ergodic result}
\label{sec:abstract}

We start by proving an abstract ergodic result, which lays the
foundations of the present work.  Recall that a Markov transition
semigroup $\CP_t$ is said to be \textit{strong Feller} at time $t$ if
$\CP_t \phi$ is continuous for every bounded measurable function
$\phi$. It is a well-known and much used fact that the strong Feller
property, combined with some irreducibility of the transition
probabilities implies the uniqueness of the invariant measure for
$\CP_t$ \cite[Th.~4.2.1]{b:DaZa96}. If $\CP_t$ is generated by a diffusion with smooth
coefficients on $\R^n$ or a finite-dimensional manifold, H\"ormander's
theorem \cite{H1}, \cite{Ho} provides us with an efficient (and sharp if the
coefficients are analytic) criterion for the strong Feller
property to hold. Unfortunately, no equivalent theorem exists if $\CP_t$
is generated by a diffusion in an infinite-dimensional space, where
the strong Feller property seems to be much ``rarer''. If the
covariance of the noise is nondegenerate (i.e.\  the diffusion is
elliptic in some sense), the strong Feller property can often be
recovered by means of the Bismut-Elworthy-Li formula \cite{Xuemei}.
The only result to our knowledge that shows the strong Feller property
for an infinite-dimensional diffusion where the covariance of the
noise does not have a dense range is given in \cite{EH3}, but it still
requires the forcing to act in a nondegenerate way on a subspace
of finite codimension.
 
\Subsec{Preliminary definitions}
Let $\X$ be a Polish (i.e.\  complete, separable, metrizable)
space. Recall that a \textit{pseudo-metric} for $\X$ is a continuous
function $d:\X^2 \to \R_+$ such that $d(x,x) = 0$ and such that the
triangle inequality is satisfied. We say that a pseudo-metric $d_1$ is
larger than $d_2$ if $d_1(x,y) \ge d_2(x,y)$ for all $(x,y) \in \X^2$.

\begin{definition}\label{def:separating}
Let $\{d_n\}_{n=0}^\infty$ be an increasing sequence of
(pseudo-)metrics on a Polish space $\X$.  If
$\lim_{n\to \infty} d_n(x,y) = 1$ for all $x \neq y$, then $\{d_n\}$ is a
\textit{totally sepa\-rating system of \/{\rm (}\/pseudo\/{\rm -)}\/metrics} for $\X$.
\end{definition}

Let us give a few representative examples.

\begin{example}
Let $\{a_n\}$ be an increasing sequence in $\R$ such that\break $\lim_{n\to
\infty} a_n= \infty$. Then, $\{d_n\}$ is a \textit{totally separating
system of \/{\rm (}\/pseudo\/{\rm -)}\/metrics} for $\X$ in the following three cases.
\begin{itemize}
\item[1.] Let $d$ be an arbitrary continuous metric on $\X$ and set
$d_n(x,y) = 1 \wedge a_n d(x,y)$.  \item[2.] Let $\X = \CC_0(\R)$ be
the space of continuous functions on $\R$ vanishing at infinity and
set $d_n(x,y) = 1 \wedge \sup_{s \in [-n,n]} a_n |x(s)-y(s)|$.
\item[3.] Let $\X = \ell^2$ and set $d_n(x,y) = 1 \wedge a_n\sum_{k
=0}^n |x_k-y_k|^2$.
\end{itemize}
\end{example}

Given a pseudo-metric $d$, we define the following seminorm on the set
of $d$-Lipschitz continuous functions from $\X$ to $\R$:
\begin{equation} \label{e:defWas}
\|\phi\|_d = \sup_{\substack{x,y \in \X\\ x \neq y}} {|\phi(x) - \phi(y)| \over d(x,y)}\;.
\end{equation}   
This in turn defines a dual seminorm on the space of
finite signed Borel measures on $\X$ with vanishing integral by
\begin{equation} \label{e:defWas2}
\norm{\nu}_d = \sup_{\|\phi\|_d = 1} \int_\X \phi(x)\,\nu(dx)\;.
\end{equation}   
Given $\mu_1$ and $\mu_2$, two positive finite Borel measures on $\X$ with equal
mass, we also denote by $\cC(\mu_1,\mu_2)$ the set of positive
measures on $\X^2$ with marginals $\mu_1$ and $\mu_2$ and we define
\begin{equation} \label{e:duality}
\|\mu_1 - \mu_2\|_d = \inf_{\mu \in \cC(\mu_1,\mu_2)} \int_{\X^2}
d(x,y)\,\mu(dx,dy)\;.
\end{equation}   
The following lemma is an easy consequence of the Monge-Kantorovich
duality; see e.g.\   \cite{Kantorovich}, \cite{Kantorovich2}, \cite{MR88f:90180}, and shows that in most cases
these two natural notions of distance can be used interchangeably.

\begin{lemma}
Let $d$ be a continuous pseudo-metric on a Polish space $\X$ and let $\mu_1$ and
$\mu_2$ be two positive measures on $\X$ with equal mass. Then{\rm ,} $\|\mu_1 - \mu_2\|_d = \norm{\mu_1 -
\mu_2}_d$.
\end{lemma}

\Proof 
This result is well-known if $(\X,d)$ is a separable metric space; see
for example \cite{MR93b:60012} for a detailed discussion on many of its
variants. If we define an equivalence relation on $\X$ by $x \sim y
\Leftrightarrow d(x,y) = 0$ and set $\X_d = \X/{\sim}$, then $d$ is
well-defined on $\X_d$ and $(\X_d,d)$ is a separable
metric space (although it may no longer be complete).  When
$\pi:\X \to \X_d$ by $\pi(x) = [x]$, the result follows from the Monge-Kantorovich 
duality in $\X_d$ and the fact that both sides of \eref{e:duality} do not change 
if the measures $\mu_i$ are replaced by $\pi^*\mu_i$.
\Endproof\vskip4pt  

Recall that the total variation norm of a finite signed measure $\mu$ on $\X$ is
given by $\|\mu\|_\TV = {1\over 2}(\mu^+(\X) + \mu^-(\X))$, where $\mu
= \mu^+ - \mu^-$ is the Jordan decomposition of $\mu$. 
The next result is crucial to the approach taken in this paper.

\begin{lemma}\label{l:limitMetrics}
Let $\{d_n\}$ be a bounded and increasing family of continuous pseudo-metrics on a Polish space $\X$ and define
$d(x,y) = \lim_{n \to \infty} d_n(x,y)$. Then{\rm ,}
  $\lim_{n \to \infty} \|\mu_1 - \mu_2\|_{d_n} = \|\mu_1 - \mu_2\|_d$ for any two positive
measures $\mu_1$ and $\mu_2$ with equal mass.
\end{lemma}

\Proof 
The limit exists since the sequence is bounded and increasing by assumption, so let us denote this limit by $L$.
It is clear from \eref{e:duality} that\break $\|\mu_1 - \mu_2\|_d \ge L$, so it remains to show the converse bound.
Let $\mu_n$ be a measure in $\cC(\mu_1,\mu_2)$ that realizes \eref{e:duality} for the distance $d_n$. (Such a measure is shown to exist in \cite{MR93b:60012}.)
The sequence $\{\mu_n\}$ is tight on $\X^2$ since its marginals are constant, and so we can extract a weakly
converging  subsequence. Denote by $\mu_\infty$ the limiting measure. For $m \geq  n$
$$
  \int_{\X^2} d_{n}(x,y)\,\mu_m(dx, dy) \leq  \int_{\X^2}
  d_{m}(x,y)\,\mu_m(dx, dy)  \leq L\;.
$$
Since $d_n$ is continuous,  the weak convergence taking $m \rightarrow \infty$ implies that
$$
\int_{\X^2} d_{n}(x,y)\,\mu_\infty(dx, dy) \le L\;,\quad \forall \, n > 0\;.
$$ 
It follows from the dominated convergence theorem that $\int_{\X^2} d(x,y)\,\mu_\infty(dx, dy)\break \le
L$, which concludes the proof.
\hfq

\begin{corollary}\label{cor:TV}
Let $\X$ be a Polish space and let $\{d_n\}$ be a totally separating
system of pseudo-metrics for $\X$. Then, $\|\mu_1-\mu_2\|_\TV =
\lim_{n \to \infty} \|\mu_1-\mu_2\|_{d_n}$ for any two positive
measures $\mu_1$ and $\mu_2$ with equal mass on $\X$.
\end{corollary}

\Proof 
It suffices to notice that $$\|\mu_1 - \mu_2\|_\TV = \inf_{\mu \in
  \cC(\mu_1,\mu_2)} \mu(\{(x,y): x \neq y\}) =  \|\mu_1 - \mu_2\|_{d}$$ with
 $d(x,y) = 1$ whenever $x\neq y$ and then to apply Lemma
  \ref{l:limitMetrics}. Observe that $d_n \to d$ by the definition of a totally separating system
  of pseudo-metrics and that Lemma~\ref{l:limitMetrics} makes
  no assumptions on the continuity of the limiting pseudo-metric $d$.
\hfq

\Subsec{Asymptotic strong Feller}
Before we define the asymptotic strong Feller property, recall that:

\begin{definition}
A Markov transition semigroup on a Polish space $\CX$ is said to be 
\textit{strong Feller} at time $t$ if
$\CP_t \phi$ is continuous for every bounded measurable function $\phi:\CX\to\R$.
\end{definition}

Note that if the transition probabilities $\CP_t(x,\cdot\,)$ are
continuous in $x$ in the total variation
topology, then $\CP_t$ is strong Feller at time $t$. 

Recall also that the support of a probability measure $\mu$, denoted by $\supp(\mu)$, is the
intersection of all closed sets of measure $1$. A useful characterization of the support of a measure
is given by

\begin{lemma}\label{lem:charge}
A point $x \in \supp(\mu)$ if and only if $\mu(U) > 0$ for every open set $U$ containing $x$. 
\end{lemma}

It is well-known that if a Markov transition semigroup $\CP_t$ is strong Feller and
$\mu_1$ and $\mu_2$ are two distinct ergodic invariant measures for $\CP_t$ (i.e.\  $\mu_1$ and $\mu_2$ are mutually singular), then
$\supp \mu_1 \cap \supp \mu_2 = \emptyset$. (This can be seen e.g.\   by the same argument
as in \cite[Prop.~4.1.1]{b:DaZa96}.)
In this section, we show that this property still holds if the strong Feller property is
replaced by the following property, where we denote by $\CU_x$ the collection of all open
sets containing $x$.

\begin{definition}\label{def:qsf}
  A Markov transition semigroup $\CP_t$ on a Polish space $\X$ is
  called {\em \qsf/} at $x$ if there exists a totally separating
  system of pseu\-do-metrics $\{d_n\}$ for $\X$ and a sequence $t_n >
  0$ such that
\begin{equation} \label{e:qsf}
  \inf_{U \in \CU_x} \limsup_{n\to\infty} \sup_{y \in U}
  \|\CP_{t_n}(x,\cdot\,) - \CP_{t_n}(y,\cdot\,)\|_{d_n} = 0\;,
\end{equation}   
It is called \qsf/ if this property holds at every $x \in \X$.
\end{definition}

\begin{remark}
If $\CB(x,\gamma)$ denotes the open ball of radius $\gamma$ centered at $x$ in some
metric defining the topology of $\X$, then it is immediate that \eref{e:qsf} is equivalent to
$$
  \lim_{\gamma \to 0} \limsup_{n\to\infty} \sup_{y \in \Ball(x,\gamma)}
  \|\CP_{t_n}(x,\cdot\,) - \CP_{t_n}(y,\cdot\,)\|_{d_n} = 0\;.
$$
\end{remark}

\begin{remark}
Notice that the definition of the asymptotic
strong Feller property allows for the possibility that $t_n = t$ for
all $n$. In this case, the transition probabilities $\CP_t (x,\cdot)$ are
continuous in the total variation topology and thus $\CP_s$ is strong
Feller at times $s \geq t$.  Conversely, strong Feller Markov semigroups
on Polish spaces are asymptotically strong Feller. To see this first
observe that if $P$ and $Q$ are two Markov operators over the same
Polish space that are strong Feller, then the product $PQ$ is a Markov
operator whose transition probabilities are continuous in the total
variation distance [DM83], [Sei02]. Hence, if $\CP_t$ is strong Feller
for some $t > 0$,  then $\CP_{2t}= \CP_t \CP_t$ is  continuous in the total
variation distance, which implies that the semigroup $\CP_t$ is asymptotically strong
Feller.  We would like to thank B. Goldys for
pointing this fact out to us.
\end{remark}

One other way of seeing the connection to the strong Feller property is to
recall that a standard criterion for $\CP_t$ to be strong Feller is given by
\cite[Lem.~7.1.5]{b:DaZa96}: 
\begin{proposition}\label{prop:classicSF}
A semigroup $\CP_t$ on a Hilbert space $\CH$ is strong Feller if{\rm ,} for all
$\phi:\CH \rightarrow \R$ with $\|\phi\|_\infty\eqdef\sup_{x \in \CH}|\phi(x)|$ and $\|\nabla \phi\|_\infty$ finite
one has
\begin{equation}
        |\nabla \CP_{t}\phi(x)| \leq C(\|x\|) \|\phi\|_\infty\;,
\end{equation}
where $C:\R_+ \rightarrow \R$ is a fixed nondecreasing function.  
\end{proposition}

The following lemma provides
a similar criterion for the asymptotic strong Feller property:
\begin{proposition} \label{l:asfLipFunctions} Let $t_n$ and $\delta_n$ be two positive sequences with
$\{t_n\}$ nondecrea\-sing and $\{\delta_n\}$ converging to zero.
A semigroup $\CP_t$ on a Hilbert space $\CH$ is asymptotically strong Feller if{\rm ,} for all
$\phi:\CH\rightarrow \R$ with $\|\phi\|_\infty$ and $\|\nabla \phi\|_\infty$ finite{\rm ,}
\begin{equation} \label{e:SF}
|\nabla \CP_{t_n}\phi(x)| \leq C(\|x\|) \bigl( \|\phi\|_\infty + \delta_n \|\nabla \phi\|_\infty\bigr) 
\end{equation}   
for all $n${\rm ,} where $C:\R_+ \rightarrow \R$ is a fixed nondecreasing function.
\end{proposition}

\Proof  
For $\eps > 0$, we define on $\CH$ the distance $$d_\eps(w_1,w_2) = 1
\wedge \eps^{-1} \|w_1-w_2\|,$$ and we denote by $\|\cdot\|_\eps$ the
corresponding seminorms on functions and on measures given by
\eref{e:defWas} and \eref{e:defWas2}. It is clear that if $\delta_n$ is
a decreasing sequence converging to $0$, $\{d_{\delta_n}\}$ is a totally
separating system of metrics for \pagebreak $\CH$.  

It follows immediately from \eref{e:SF} that for every
Fr\'echet differentiable function $\phi$ from $\CH$ to $\R$ with
$\|\phi\|_\eps \le 1$,
\begin{equation} \label{e:mainbound2}
\int_\CH \phi(w)\,\bigl(\CP_{t_n}(w_1,dw) - \CP_{t_n}(w_2,dw)\bigr)\le \|w_1-w_2\|
C(\|w_1\|\vee\|w_2\|)\Bigl(1 + \frac{ \delta_n}{\eps}\Bigr)\;.
\end{equation}   

Now take a Lipschitz continuous function $\phi$ with $\|\phi\|_\eps \le 1$.
By applying to $\phi$ the semigroup at time $1/m$ corresponding to a linear
strong Feller diffusion in $\CH$, one obtains (\cite{CerrRD}, \cite{b:DaZa96}) a sequence $\phi_m$ of
Fr\'echet differentiable approximations $\phi_m$ with $\|\phi_m\|_\eps \le 1$ and
such that $\phi_m \to \phi$ pointwise. Therefore, by the dominated convergence
theorem, \eref{e:mainbound2} holds for Lipschitz continuous functions $\phi$ 
and so 
$$
\|\CP_{t_n}(w_1,\cdot\,) - \CP_{t_n}(w_2,\cdot\,)\|_{\eps}\le \|w_1-w_2\|
C(\|w_1\|\vee\|w_2\|)\Bigl(1 + \frac{\delta_n}\eps\Bigr)\;.
$$   
Choosing $\eps=a_n=\sqrt{\delta_n}$, we obtain
$$
\|\CP_{t_n}(w_1,\cdot\,) - \CP_{t_n}(w_2,\cdot\,)\|_{a_n}\le \|w_1-w_2\|
C(\|w_1\|\vee\|w_2\|)\bigl(1 + a_n\bigr)\;,
$$  
which in turn implies that $\CP_t$
is \qsf/ since $a_n \rightarrow 0$.
\hfq

\begin{example}
Consider the SDE
$$
  dx = -x\,dt + dW(t)\;,\qquad dy = -y\,dt\;.
$$
Then, the corresponding Markov semigroup $\CP_t$ on $\R^2$ is not
strong Feller, but it is \qsf/. To see that $\CP_t$ is not strong
Feller, let $\varphi(x,y)=\sgn(y)$ and observe that $\CP_t \varphi =
\varphi$ for all $t \in [0,\infty)$. Since $\varphi$ is bounded but
not continuous, the system is not strong Feller. To see that the
system is \qsf/ observe that for any differentiable $\varphi:\R^2
\rightarrow \R$ and any direction $\xi \in \R^2$ with $\|\xi\|=1$,  
\begin{eqnarray*}
  \big|(\nabla\CP_t\varphi)(x_0,y_0) \cdot \xi\big| &=& \big|
  \EE_{(x_0,y_0)}(\nabla \varphi)(x_t,y_t)\cdot (u_t,v_t)\big| \\
  &\leq&  \|\nabla \varphi\|_\infty \EE \big|(u_t,v_t)\big| \leq
  \|\nabla \varphi\|_\infty e^{-t}\;,
\end{eqnarray*}
where $(u_t,v_t)$ is the linearized flow starting from $\xi$. In other
words $(u_0,v_0)=\xi$, $d u= -udt$, and $d v = -vdt$. This  is
a particularly simple example  because the flow is globally contractive. 
\end{example}

\begin{example}
  Now consider the SDE
$$
  dx = (x-x^3)\,dt + dW(t)\;,\qquad dy = -y\,dt\;.
$$ 
Again the function $\varphi(x,y)=\sgn(y)$ is invariant under $\CP_t$
implying that the system is not strong Feller. It is however
not globally contractive.  As in the
previous example, let $\xi=(\xi_1,\xi_2) \in \R^2$   and $\|\xi\|=1$
and now let $(u_t,v_t)$ denote the linearizion of this equation with
$(u_0,v_0)=\xi$. Let $\CP_t^x$ denote the Markov transition semigroup
of the $x_t$ process. It is a classical fact that for such a uniformly
elliptic diffusion with a unique invariant measure one has $|\d_x\CP_t^x\varphi(x,y)| \leq
C(|x|)\|\varphi\|_\infty$ for some nondecreasing function $C$ and
all \pagebreak $t\geq 1$.
Hence differentiating with respect to both initial conditions produces
\begin{eqnarray*}
\big| (\nabla\CP_t\varphi)(x_0,y_0) \cdot \xi\big| &= \big|
  \big(\d_x\CP_t^x\varphi(x,y) \xi_1 \big) +
  \EE\big((\d_y\varphi)(x_t,\tilde y_t)v_t\big)\big|\\
  &\leq C(|x|)\|\varphi\|_\infty + \EE\big| v_t\big| \|\nabla
  \varphi\|_\infty \\
    &\leq \big(C(|x|)+1\big)\big(  \|\varphi\|_\infty  +  e^{-t} \|\nabla
  \varphi\|_\infty \big)
\end{eqnarray*}
for $t \geq 1$ which implies that the system is \qsf/. 
\end{example}

\begin{example}
  In infinite dimensions, even a seemingly nondegenerate diffusion
  can suffer from a similar problem. Consider the following infinite
  dimensional Ornstein-Uhlenbeck process $u(x,t)=\sum \hat u(k,t)
  \exp(ikx)$ written in terms of its complex Fourier coefficients. We
  take $x \in \mathbb{T}=[-\pi,\pi]$, $k \in \Z$ and 
\begin{equation} \label{e:OU}
    d \hat u(k,t) = -(1+|k|^2) \hat u(k,t)\,dt + \exp(-|k|^3)\, d\beta_k(t)\;,
\end{equation}   
where the $\beta_k$ are independent standard complex Brownian motions.
The Markov transition densities $\CP_t(x,\cdot)$ and $\CP_t(y,\cdot)$
are singular for all finite times if $x-y$ is not sufficiently smooth.
This implies that the diffusion \eref{e:OU} in $\CH =
\L^2([-\pi,\pi])$ is \textit{not} strong Feller since by Lemma 7.2.1
of \cite{b:DaZa96} the strong Feller property is equivalent to
$\CP_t(y,\cdot)$ being equivalent to $\CP_t(x,\cdot)$ for all $x$ and
$y$. Another equivalent characterization of the strong Feller property
is that the $\mbox{image}(S_t) \subset \mbox{image}(Q_t)$ where $S_t$
is the linear semigroup generated by the deterministic part for the
equation defined by $(S_tu)(k)=e^{-(1+|k|^2)t}u(k,0)$ and $Q_t=\int_0^t
S_r G S_r^*dr$ where $G$ is the covariance operator of the noise
defined by $(Gu)(k)= \exp(-2|k|^3)u(k)$. This captures the fact
that the mean, controlled by $S_t$, is moving towards zero too slowly
relative to the decay of the noise's covariance structure. 
However, one can easily check that the example is \qsf/ since the entire flow
is pathwise contractive as in the first example.
\end{example}

The classical strong Feller property captures well the smoothing due
to the random effects. When combined with
irreducibility in the same topology, it implies that the
transition densities starting from different points are mutually
absolutely continuous. As the examples show, this is often not true in
infinite dimensions. We see that
the asymptotic strong Feller property better incorporates the smoothing due to the pathwise
contraction of the dynamics. Comparing Proposition \ref{prop:classicSF} with
Proposition \ref{l:asfLipFunctions}, one sees that the second term in
Proposition \ref{l:asfLipFunctions} allows one to capture the
progressive smoothing in time from the pathwise dynamics. This becomes
even clearer when one examines the proofs of Proposition
\ref{l:theWholeEnchilada} and Proposition \ref{effectivelyEllipticASF}
later in the text. There one sees that the first term comes from
shifting a derivative from the test function to the Wiener measure and
the second is controlled using  in an essential way the
contraction due to the spatial Laplacian.

The usefulness of the \qsfp/ is seen in the following theorem
and its accompanying corollary which are the main results of this
section.
\begin{theorem}\label{theo:quasiStrongFeller}
Let $\CP_t$ be a Markov semigroup on a Polish space $\X$ and let $\mu$ and $\nu$ be two
distinct ergodic invariant probability measures for $\CP_t$. If $\CP_t$ is \qsf/ at $x${\rm ,} then $x \not\in
\supp \mu \cap \supp \nu$.
\end{theorem}

\Proof 
  By \cor{cor:TV}, the proof of this result is a simple rewriting
  of the proof of the corresponding result for strong Feller
  semigroups.
  
   For every measurable set $A$, every $t>0$, and every pseudo-metric
  $d$ on $\X$ with $d \leq 1$, the triangle inequality for $\|\cdot\|_d$ implies
\begin{equation} \label{e:ineq}
  \|\mu - \nu\|_d \le 1- \min\{\mu(A), \nu(A)\} \Bigl(1-\max_{y,z \in
    A}\|\CP_t(z,\cdot) - \CP_t(y,\cdot)\|_d\Bigr)\;.
\end{equation}   
To see this, set $\alpha = \min\{ \mu(A),\nu(A)\}$. If $\alpha =0$
there is nothing to prove so assume $\alpha >0$.  Clearly there exist
probability measures $\bar \nu$, $\bar \mu$, $\nu_A$, and $\mu_A$ such that
$\nu_A(A) = \mu_A(A) = 1$ and such that 
$\mu=(1-\alpha) \bar\mu+ \alpha \mu_A$ and $\nu=(1-\alpha)\bar \nu +
\alpha \nu_A$. From  the invariance of the measures $\mu$ and $\nu$ and the triangle
inequality this  implies
\begin{align*}
  \|\mu - \nu\|_d &= \|\CP_t\mu - \CP_t\nu\|_d \leq
  (1-\alpha)\|\CP_t\bar\mu - \CP_t\bar\nu\|_d + \alpha \|\CP_t\mu_A -
  \CP_t\nu_A\|_d\\
  &\leq (1-\alpha) + \alpha \int_A \int_A \|\CP_t(z,\cdot)  -
  \CP_t(y,\cdot)\|_d \mu_A(dz) \nu_A(dy)\\
  &\leq  1- \alpha \Bigl( 1-\max_{y,z \in A} \| \CP_t(z,\cdot) -
  \CP_t(y,\cdot)  \|_d   \Bigr)\; .
\end{align*}
 
Continuing with the proof of the corollary, we see that, by the definition of the
\qsfp/, there exist constants $N>0$, a sequence of totally separating
pseudo-metrics $\{d_n\}$, and an
open set $U$ containing $x$ such that $\|\CP_{t_n}(z,\cdot) -
\CP_{t_n}(y,\cdot)\|_{d_n} \le {1/2}$ for every $n > N$ and every $y,z
\in U$. (Note that by the definition of totally separating
pseudo-metrics $d_n \leq 1$.)  

Assume by contradiction that $x \in \supp \mu \cap \supp\nu$ and
therefore that $\alpha=\min(\mu(U),\nu(U))>0$.  Taking $A = U$,
$d=d_n$, and $t = t_n$ in \eref{e:ineq}, we then get $\|\mu -
\nu\|_{d_n} \le 1- {\alpha\over 2}$ for every $n > N$, and therefore
$\|\mu - \nu\|_\TV \le 1-{\alpha\over 2}$ by \cor{cor:TV}, thus
leading to a contradiction.
\Endproof\vskip4pt  

As an immediate corollary, we have

\begin{corollary}\label{cor:main}
  If $\CP_t$ is an \qsf/ Markov semigroup and there exists a point $x$ such
  that $x \in \supp \mu$ for every invariant probability measure $\mu$
  of $\CP_t${\rm ,} then there exists at most one invariant probability
  measure for $\CP_t$.
\end{corollary}

\section{Applications to the stochastic 2D Navier-Stokes equations}
\label{sec:SNS}

To state the general ergodic result for the two-dimensional
Navier-Stokes equations, we begin by looking at the 
algebraic structure of the Navier-Stokes nonlinearity   in
Fourier space.

Remember that $\CZ_0$ as given in \eref{e:defZ0} denotes the set of
forced Fourier modes for \eref{e:vort}. In view of
Equation~\ref{e:vortFourier}, it is natural to consider the set
$\tilde \CZ_\infty$, defined as the smallest subset of $\Z^2$
containing $\CZ_0$ and satisfying that for every $\ell,j \in \tilde
\CZ_\infty$ such that $\scal{\ell^\perp, j} \not =0$ and $|j|\not
=|\ell|$, one has $j+\ell \in \tilde \CZ_\infty$
(see \cite{b:EMattingly00}). Denote by $\tilde \CH$ the closed subspace of
$\CH$ spanned by the Fourier basis vectors corresponding to elements
of $\tilde \CZ_\infty$.  Then, $\tilde \CH$ is invariant under the
flow defined by \eref{e:vort}.

Since we would like to make use of the existing results, we recall
the sequence of subsets $\CZ_n$ of $\Z^2$ defined recursively
in \cite{b:MattinglyPardoux03Pre} by
$$
  \mathcal{Z}_n =  \Big\{ \ell+j \,\Big|\,j \in
  \mathcal{Z}_0, \ell\in \mathcal{Z}_{n-1}\, \text{with}\,
  \scal{\ell^\perp, j} \not =0, |j|\not =|\ell| \Big\}\; ,
$$   
as well as
$ \mathcal{Z}_\infty=\bigcup_{n=1}^\infty \mathcal{Z}_n$.
The two sets $\CZ_\infty$ and $\tilde \CZ_\infty$ are the same even
though from the definitions we only see $\CZ_\infty \subset \tilde
\CZ_\infty$. The other inclusion follows from the characterization of
$\CZ_\infty$ given in \prop{prop:gen} below.

The following theorem is the principal result of this article.

\begin{theorem}\label{theo:main}
The transition semigroup on $\tilde \CH$ generated by the solutions to \eref{e:vort} 
is \qsf/.
\end{theorem}

An almost immediate corollary of \theo{theo:main} is

\begin{corollary}\label{c:ergodic}
There exists exactly one invariant probability measure for \eref{e:vort} restricted to $\tilde \CH$.
\end{corollary}

{\it Proof of Corollary {\rm \ref{c:ergodic}.}}
  The existence of an invariant probability measure $\mu$ for
  \eref{e:vort} is a standard result \cite{b:Fl94}, \cite{b:DaZa96}, \cite{b:ChowKhasminskii98}. 
 By \cor{cor:main} it suffices to show that the
  support of every invariant measure contains the element $0$.  Applying It\^o's
  formula to $\|w\|^2$ yields for every invariant measure $\mu$ the
  {\it a~priori} bound
$$
 \int_\CH \|w\|^2\,\mu(dw) \le \frac{C\mathcal{E}_0}{\nu}\;.
$$  
(See \cite[Lemma B.1]{b:EMattinglySinai00}.) Therefore, denoting by
$\Ball(\orho)$ the ball of radius $\orho$ centered at $0$, we have $\tilde C$
such that $\mu\bigl(\Ball(\tilde C)\bigr) > {1\over 2}$ for every
invariant measure~$\mu$. On the other hand,
\cite[Lemma~3.1]{b:EMattingly00} shows that, for every $\gamma > 0$
there exists a time $T_\gamma$ such that
$$
\inf_{w \in \Ball(\tilde C)} \CP_{T_\gamma}\bigl(w,\Ball(\gamma)\bigr) > 0\;.
$$   
(Note, though \cite[Lemma~3.1]{b:EMattingly00} was about Galerkin
approximations, inspection of the proof reveals that it holds equally
for the full solution.)
Therefore, $\mu(\Ball(\gamma)) > 0$ for every $\gamma >0$ and every
invariant measure $\mu$, which implies that $0 \in \supp(\mu)$ by
\lem{lem:charge}. 
\Endproof\vskip4pt

The crucial ingredient in the proof of \theo{theo:main} is the
following result:

\begin{proposition} \label{l:theWholeEnchilada} 
For every $\eta > 0${\rm ,} there exist 
constants $C, \delta > 0$ such
that for every Fr\ee chet differentiable function $\phi$ from $\tilde \CH$
to $\R$ one has the bound 
\begin{equation} \label{e:mainbound}
\|\nabla\CP_n \phi(w)\|\le C\exp(\eta\|w\|^2)\bigl(\|\phi\|_\infty
      + \|\nabla\phi\|_\infty e^{-\delta n}\bigr)\;,
\end{equation}   
for every $w \in \tilde \CH$ and $n \in \N$.
\end{proposition}

The proof of Proposition~\ref{l:theWholeEnchilada} is the content of
Section~\ref{sec:proof} below. \theo{theo:main} then follows from this proposition
and from Proposition
\ref{l:asfLipFunctions} with the choices $t_n = n$ and
$\delta_n=e^{-\delta n}$.
 Before we turn to the proof of Proposition~\ref{l:theWholeEnchilada},  we  characterize
$\mathcal{Z}_\infty$ and give an
informal introduction to Malliavin calculus adapted to our framework, followed by a 
brief discussion on how it relates to the
strong Feller property.

\Subsec{The structure of $\mathcal{Z}_\infty$}
\label{sec:structureZinfty}
In this section, we give a complete characterization of the set $\CZ_\infty$. 
We start by defining $\scal{\CZ_0}$ as the subset of $\Z^2\setminus \{(0,0)\}$ generated 
by integer linear  combinations of elements of $\CZ_0$. With this notation, we have

\begin{proposition}\label{prop:gen}
If there exist $a_1, a_2 \in \CZ_0$ such that $|a_1| \neq |a_2|$ and such that $a_1$ and
$a_2$ are not collinear{\rm ,} then $\CZ_\infty = \scal{\CZ_0}$. Otherwise, $\CZ_\infty = \CZ_0$.
In either case{\rm ,} one always has that $\CZ_\infty = \tilde \CZ_\infty$.
\end{proposition}

This also allows us to characterize the main case of interest:

\begin{corollary}\label{cor:gen}
  One has $\mathcal{Z}_\infty = \Z^2 \setminus \{(0,0)\}$ if and only
  if the following holds\/{\rm :}\/
\begin{itemize}
\item[{\rm 1.}] Integer linear combinations of elements of $\CZ_0$ generate
  $\Z^2$.
\item[{\rm 2.}] There exist at least two elements in $\CZ_0$ with nonequal
  Euclidean norm.
\end{itemize}
\end{corollary}

\demo{Proof of \prop{prop:gen}}
It is clear from the definitions that if the elements of $\CZ_0$ are
all collinear or of the same Euclidean length, one has $\CZ_\infty = \CZ_0 = \tilde \CZ_\infty$. 
In the rest of the proof, we assume that there exist two elements $a_1$ and $a_2$ 
of $\CZ_0$ that are neither collinear nor of the same length and we show that one has
$\CZ_\infty = \scal{\CZ_0}$. Since it follows from the definitions that $\CZ_\infty \subset \tilde \CZ_\infty \subset \scal{\CZ_0}$, this shows that $\CZ_\infty = \tilde \CZ_\infty$.

Note that the set $\CZ_\infty$ consists exactly of those points in $\Z^2$ that can be
reached by a walk starting from the origin with steps drawn in $\CZ_0$ and which does not
contain any of the following ``forbidden steps'':
\begin{definition}
  A step with increment $\ell \in \CZ_0$ starting from $j \in \Z^2$ is
  \textit{forbidden} if either $|j| = |\ell|$ or $j$ and $\ell$ are
  collinear.
\end{definition}

Our first aim is to show that there exists $R>0$ such that
$\CZ_\infty$ contains every element of $\scal{\CZ_0}$ with Euclidean 
norm larger than $R$. In order to achieve this, we start with a few very simple
observations.

\begin{lemma}\label{lem:trivial}
  For every $R_0 > 0$, there exists $R_1 > 0$ such that every $j\in
  \scal{\CZ_0}$ with $|j| \le R_0$ can be reached from the origin by a path
  with steps in $\CZ_0$ \/{\rm (}\/some steps may be forbidden\/{\rm )}\/ which never
  exits the ball of radius $R_1$.  
\end{lemma}

\begin{lemma}\label{lem:grid}
  There exists $L>0$ such that the set $\CZ_\infty$ contains all
  elements of the form $n_1 a_1 + n_2 a_2$ with $n_1$ and $n_2$ in $\Z
  \setminus [-L,L]$.
\end{lemma}

\Proof 
  We may assume without loss of generality that $|a_1| \!>\! |a_2|$ and that
  $\scal{a_1,a_2} > 0$.  Choose $L$ such that $L \scal{a_1,a_2}\! \ge\!
  |a_1|^2$.  By the symmetry of $\CZ_0$, we can replace $(a_1,a_2)$ by
  $(-a_1,-a_2)$, so that we can  assume without loss of generality that
  $n_2 > 0$. We then take first one step in the direction $a_1$ starting from the origin,
  followed by $n_2$ steps in the direction $a_2$. Note that the
  assumptions we made on $a_1$, $a_2$, and $n_2$ ensure that none of
  these steps is forbidden. From there, the condition $n_2 > L$
  ensures that we can take as many steps as we want into either the
  direction $a_1$ or the direction $-a_1$ without any of them being
  forbidden.
\Endproof\vskip4pt  

Denote by $Z$ the set of elements of the form $n_1a_1 + n_2 a_2$ considered in
Lem\-ma~\ref{lem:grid}.  It is clear that there exists $R_0 > 0$ such that
every element in $\scal{\CZ_0}$ is at distance less than $R_0$ of an element
of $Z$. Given this value $R_0$, we now fix $R_1$ as given from
\lem{lem:trivial}.  Let us define the set
$$
  A = \Z^2 \cap \bigl(\{\alpha j \,|\, \alpha \in \R\,,\; j \in
  \CZ_0\} \cup \{k \,|\, \exists j \in \CZ_0 \,\text{with}\, |j| =
  |k|\}\bigr)\;,
$$
which has the property that there is no forbidden step starting from
$\Z^2 \setminus A$.  Define furthermore
$$
B = \{j \in \scal{\CZ_0}\,|\, \inf_{k \in A} |k-j| > R_1\}\;.
$$

\begin{figure}
\centerline{\epsfig{file=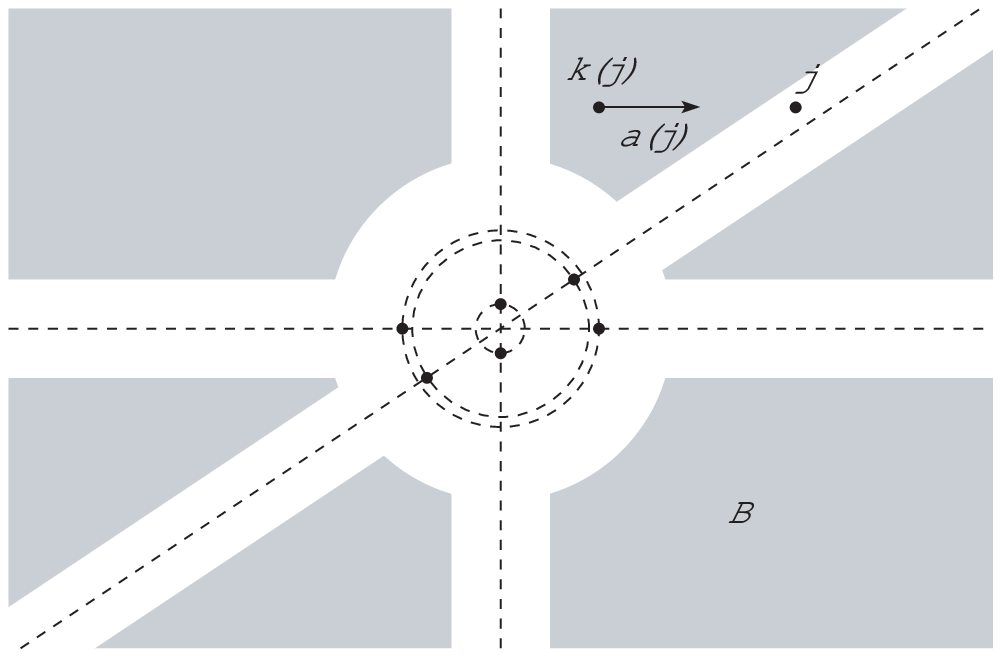}}
  \caption{Construction from the proof of \prop{prop:gen}.}
  \label{fig:sets}
\end{figure}

By \lem{lem:trivial} and the definition of $B$, every element of $B$ 
can be reached by a path from $Z$ containing no forbidden steps, therefore $B \subset \CZ_\infty$.
On the other hand, it is easy to see that there exists $R > 0$ such that for every
element of $j \in \scal{\CZ_0} \setminus B$ with $|j| > R$, there exists an
element $a(j) \in \CZ_0$ and an element $k(j) \in B$ such that $j$ can
be reached from $k(j)$ with a finite number of steps in the direction
$a(j)$. Furthermore, if $R$ is chosen sufficiently large, none of these
steps crosses $A$, and therefore none of them is forbidden.
We have thus shown that there exists
$R>0$ such that $\CZ_\infty$ contains $\{j\in \scal{\CZ_0}\,|\, |j|^2 \ge
R\}$. 

In order to help in  visualizing this construction, Figure~\ref{fig:sets} shows the
typical shapes of the sets $A$ (dashed lines) and $B$ (gray area), as
well as possible choices of $a(j)$ and $k(j)$, given $j$. (The black
dots on the intersections of the circles and the lines making up $A$
depict the elements of $\CZ_0$.)

We can (and will from now on) assume that $R$ is an integer. The
last step in the proof of \prop{prop:gen} is

\begin{lemma}\label{lem:recursion}
Assume that there exists an integer $R>1$ such that $\CZ_\infty$
contains $\{j\in \scal{\CZ_0}\,|\, |j|^2 \ge R\}$.  Then $\CZ_\infty$ also
contains $\{j\in \scal{\CZ_0}\,|\, |j|^2 \ge R-1\}$.
\end{lemma}

\Proof 
Assume that the set $\{j\in \scal{\CZ_0}\,|\, |j|^2 = R-1\}$ is nonempty
and choose an element $j$ from this set. Since $\CZ_0$ contains at least two
elements that are not collinear, we can choose $k \in \CZ_0$ such that $k$ is
not collinear to $j$.  Since $\CZ_0$ is closed under the operation $k
\mapsto -k$, we can assume that $\scal{j,k} \ge 0$. Consequently, one has 
$|j+k|^2 \ge R$, and so $j+k \in
\CZ_\infty$ by assumption. The same argument shows that $|j+k|^2
\ge |k|^2 + 1$, so the step $-k$ starting from $j+k$ is not forbidden and
therefore $k\in\CZ_\infty$.
\Endproof\vskip4pt  

This shows that $\CZ_\infty=\scal{\CZ_0}$ and therefore completes the proof
of \prop{prop:gen}. 
\hfq

\Subsec{Malliavin calculus and the Navier-Stokes equations}
In this section, we give a brief introduction to some elements of Malliavin
calculus applied to equation \eref{e:vort} to help orient the reader
and fix notation. We refer to \cite{b:MattinglyPardoux03Pre} for a longer
introduction in the setting of equation \eref{e:vort} and to
\cite{MR96k:60130}, \cite{b:Bell87} for a more general introduction.

Recall from Section \ref{sec:basicSettingMainResults}, that  $\Phi_t\colon C([0,t];\R^m) \times
\CH \rightarrow \CH$ was the map such that  $w_t=\Phi_t(W, w_0)$  for
initial condition $w_0$ and noise realization $W$.
Given a $v \in \L^2_{\rm loc}(\R_+, \R^m)$, the Malliavin derivative of
the $\CH$-valued random variable $w_t$ in the direction $v$, denoted
$\cD^v w_t$, is defined by
$$
 \cD^v w_t =\lim_{\eps \rightarrow 0} \frac{\Phi_t(W+\eps V, w_0)
 -\Phi_t(W,w_0)}{\eps}\;,
$$   
where the limit holds almost surely with respect to the Wiener measure
and where we set $V(t)=\int_0^t v(s)\,ds$.  Note that we allow $v$ to be 
random and possibly nonadapted
to the filtration generated by the increments of $W$.

Defining the symmetrized nonlinearity $\tilde B(w,v) = B(\KK w,v) +
B(\KK v,w)$, we use the notation $\J_{s,t}$ with $s\le t$ for the
derivative flow between times $s$ and~$t$, i.e.\  for every $\uu \in \CH$,
$\J_{s,t}\uu$ is the solution of
\begin{equation}\label{eq:JLinear}
\d_t \J_{s,t}\uu = \nu \Delta \J_{s,t}\uu + \tilde B(w_t, \J_{s,t}\uu),\quad
t>s\;,\quad \J_{s,s}\uu = \uu\;.
\end{equation}
Note that we have the important cocycle property $\J_{s,t} =
\J_{r,t}\J_{s,r}$ for $r\in [s,t]$.

Observe that $\cD^v w_t= A_{0,t}v$ where the random operator
$A_{s,t}:\L^2([s,t],\R^m)\break\to\CH$ is given by
$$
  A_{s,t}v=\int_s^t \J_{r,t}Q v(r)\,dr \;.
$$
To summarize, $\J_{0,t}\xi$ is the effect on $w_t$ of an infinitesimal
perturbation of the initial condition in the direction $\xi$ and
$A_{0,t}v$ is the effect on $w_t$ of an infinitesimal perturbation of
the Wiener process in the direction of $V(s)=\int_0^s v(r)\,dr$.

Two fundamental facts we will use from Malliavin calculus
are embodied in the following equalities. The first amounts to the
chain rule, the second is integration by parts. For a smooth function
$\phi:\CH \rightarrow \R$ and a (sufficiently regular) process $v$,
\begin{equation}\label{e:integrationByParts}
 \EE \scal{(\nabla\phi)(w_t), \cD^vw_t}=\EE \Big(\cD^v
 \bigl(\phi(w_t)\bigr)\Big)=\EE\Big( \phi(w_t) \int_0^t \scal{v(s), dW_s} \Big)\;.
\end{equation}
The stochastic integral appearing in this expression is an It\^o
integral if the process $v$ is adapted to the filtration $\CF_t$ generated
by the increments of $W$ and a
Skorokhod integral otherwise.

We also need the adjoint $A_{s,t}^*:\CH~\to~\L^2([s,t],\R^m)$ defined
by the duality relation $\scal{A_{s,t}^*\xi,v} = \scal{\xi, A_{s,t}v}$, where the
first scalar product is in $\L^2([s,t],\R^m)$ and the second one is in $\CH$.
Note that one has $(A^*_{s,t} \xi)(r)=Q^* \J_{r,t}^* \xi$, where $J_{r,t}^*$ is
the adjoint in $\CH$ of $J_{r,t}$.

One of the fundamental objects in the study of hypoelliptic diffusions
is the Malliavin matrix $M_{s,t}\eqdef A_{s,t}^{} A_{s,t}^*$. A
glimpse of its importance can be seen from the following. For $\xi \in
\CH$,  
\begin{equation*}
   \scal{M_{0,t}\xi,\xi}= \sum_{i=1}^m\int_0^t \scal{\J_{s,t} Qe_i,\xi}^2\, ds
   \ .
\end{equation*}
Hence the quadratic form $\scal{M_{0,t}\xi,\xi}$ is zero for a direction
$\xi$ only if no variation whatsoever in the Wiener process at times $s
\leq t$ could cause a variation in $w_t$ with a nonzero component in
the direction $\xi$.

We also recall that the second derivative $\K_{s,t}$ of the flow is the
bilinear map solving
\begin{eqnarray*}
\d_t \K_{s,t}(\uu,\uu') &=& \nu \Delta \K_{s,t}(\uu,\uu') + \tilde B(w_t,
\K_{s,t}(\uu,\uu')) + \tilde B(\J_{s,t}\uu', \J_{s,t}\uu)\;,\\ \K_{s,s}(\uu,\uu')
&=& 0\;.
\end{eqnarray*}
It follows from the variation-of-constants formula that
$\K_{s,t}(\uu,\uu')$ is given by
\begin{equation} \label{e:defK}
\K_{s,t}(\uu,\uu') = \int_s^t \J_{r,t}\tilde B(\J_{s,r}\uu',
\J_{s,r}\uu)\,dr\;.
\end{equation}   


\Subsec{Motivating discussion}
\label{sec:wishfulThinking}
It is instructive to proceed formally pretending that $M_{0,t}$ is
invertible as an operator on $\CH$. This is probably not true for
the problem considered here and we will certainly not attempt to prove
it in this article, but the proof presented in Section~\ref{sec:proof} is a
modification of the argument in the invertible case and hence it is
instructive to start there.

Set  $\uu_t=\J_{0,t}\xi$; now  $\xi_t$ can be interpreted as the
perturbation of $w_t$ caused by a perturbation $\xi$ in the initial
condition of $w_t$.  Our goal is to find an infinitesimal variation in
the Wiener path $W$ over the interval $[0,t]$ which produces the same
perturbation at time $t$ as the shift in the initial condition. We
want to choose the variation which will change the value of the
density the least. In other words, we choose the path with the least
action with respect to the metric induced by the inverse of the Malliavin matrix. The
least squares solution to this variational problem is easily seen to
be, at least formally, $v=A_{0,t}^* M_{0,t}^{-1} \uu_t$ where $v \in
\L^2([0,t],\R^m)$.  Observe that $\cD^v w_t = A_{0,t}v=\J_{0,t}\xi$.
Considering the derivative with respect to the initial condition $w$
of the Markov semigroup $\CP_t$ acting on a smooth function $\phi$,
we obtain
\begin{align}\label{eq:intbyParts}
  \scal{\nabla \CP_t\phi(w),\xi} &=
  \EE_{w}\bigl((\nabla\phi)(w_t)\J_{0,t}\xi\bigr)= \EE_{w}\bigl((\nabla\phi)(w_t)\cD^v
  w_t\bigr)\\&=\EE_{w}\bigg( \phi(w_t) \int_0^t v(s) dW_s \bigg)\leq
  \|\phi\|_\infty \EE_{w} \left| \int_0^t v(s) dW_s \right|\;,\notag
\end{align}
where the penultimate estimate follows from the integration by parts formula
\eref{e:integrationByParts}. Since the last term in the
chain of implications holds for functions which are simply bounded and
measurable, the estimate extends by approximation to that class of
$\phi$. Furthermore since the constant $\EE_{w} \bigl| \int_0^t v(s)
  dW_s \bigr|$ is independent of $\phi$, if one can show it is finite
and bounded independently of $\xi \in \tilde \CH$ with $\|\xi\| = 1$, we have 
proved that $\|\nabla \CP_t\phi\|$ is bounded and thus that $\CP_t$ is 
strong Feller in the 
topology of $\tilde \CH$. Ergodicity then follows from this statement by means of
\cor{cor:main}. In particular, the estimate in \eref{e:mainbound}
would hold.

In   slightly different language, since $v$ is the infinitesimal shift
in the Wiener path equivalent to the infinitesimal variation in the
initial condition~$\xi$, we have, via the Cameron-Martin
theorem, the infinitesimal change in the Radon-Nikodym derivative of
the ``shifted'' measure with respect to the original Wiener measure.
This is not trivial since in order to compute the shift~$v$, one uses
information on $\{w_s\}_{s \in [0,t]}$, so that it is in general not
adapted to the Wiener process $W_s$. This nonadaptedness can be
overcome as Section~\ref{sec:costControl} demonstrates. However the assumption
in the above calculation that $M_{0,t}$ is invertible is more serious.
We will overcome this by using the ideas and understanding which begin
in \cite{b:Mattingly98b}, \cite{b:Mattingly98}, \cite{b:EMattinglySinai00}, \cite{b:KuksinShirikyan00},
\cite{b:BricmontKupiainenLefevere01}.

The difficulty in inverting $M_{0,t}$ partly lies in our incomplete
understanding of the natural space in which \eqref{e:vort} lives.  The
knowledge needed to identify on what domain $M_{0,t}$ can be inverted
seems equivalent to identifying the correct reference measure against
which to write the transition densities. By ``reference measure,'' we
mean a replacement for the role of Lebesgue measure from finite
dimensional diffusion theory. This is a very difficult proposition. An
alternative was given in the papers
\cite{b:Mattingly98b},
\cite{b:Mattingly98},
\cite{b:KuksinShirikyan00},
\cite{b:EMattinglySinai00},
\cite{b:BricmontKupiainenLefevere01},
\cite{MatNS}, \cite{b:BricmontKupiainenLefevere02}, \cite{HExp02}, \cite{b:MasmoudiYoung02}. The idea was
to use the pathwise contractive properties of the flow at small scales due to the presence of the spatial Laplacian.
Roughly speaking, the system has finitely many unstable directions and infinitely many stable directions. One can
then use the noise to steer the unstable directions together and let the dynamics cause the stable directions to contract.
This requires the small scales to be enslaved to the large scales in some sense. A stochastic version of such a
determining modes statement (cf.\  \cite{b:FoiasProdi67}) was developed in \cite{b:Mattingly98b}.  Such an
approach to prove ergodicity requires looking at the entire future to $+\infty$ (or equivalently the entire past) as the
stable dynamics only brings solutions together asymptotically.  In the first
papers in the continuous time setting
(\cite{b:EMattinglySinai00}, \cite{MatNS}, \cite{b:BricmontKupiainenLefevere02}),
Girsanov's theorem was used to bring the unstable directions together
completely; \cite{HExp02} demonstrated the effectiveness of only
steering all of the modes together asymptotically. Since all of these
techniques used Girsanov's theorem, they required that all of the
unstable directions be directly forced. This is a type of partial
ellipticity assumption, which we will refer to as ``effective
ellipticity.'' The main achievement of this text is to remove this
restriction. We also make another innovation which simplifies the
argument considerably. We work infinitesimally, employing the
linearization of the solution rather than looking at solutions
starting from two different starting points.

\Subsec{Preliminary calculations and discussion}
\label{sec:preliminaryCalc}
Throughout this and the following sections we fix once and for all the initial
condition $w_0 \in \tilde \CH$ for \eref{e:vort} and denote by $w_t$ the
stochastic process solving \eref{e:vort} with initial condition
$w_0$. By $\EE$ we mean the expectation starting from this
initial condition unless otherwise indicated. Recall also the notation
$\CE_0 = \tr QQ^*=\sum |q_k|^2$. 
The following lemma provides  the auxiliary estimates which will be used to control various
terms during the proof of \prop{l:theWholeEnchilada}.
 
\begin{lemma}\label{lem:mainbounds}
The solution of the {\rm 2D} Navier-Stokes equations in the vorticity
formulation \eref{e:vort} satisfies the following bounds\/{\rm :}\/
\begin{itemize}
\ritem{1.}  There exist positive constants $C$ and
  $\eta_0${\rm ,} depending only on $Q$ and $\nu${\rm ,} such that \end{itemize}\vglue-30pt
\begin{multline}\label{e:boundInt}
\expect \exp\Bigl( \eta \sup_{t \geq s} \Bigl(\|w_t\|^2+ \nu  \int_s^t \|w_r\|_1^2\,dr -\mathcal{E}_0 (t-s)\Bigr)\Bigr)
\\\le C \exp \bigl(\eta e^{-\nu s}\|w_0\|^2\bigr)\;,
\end{multline} \begin{itemize}\item[]
for every $s \ge 0$ and for every $\eta \le \eta_0$. Here and in the sequel{\rm ,}   the notation $\|w\|_1 = \|\nabla w\|$
is used.
\ritem{2.}  There exist constants $\eta_1, a, \gamma > 0${\rm ,}
  depending only on $\mathcal{E}_0$ and $\nu${\rm ,} such that \end{itemize} \vglue-20pt
\begin{equation} \label{e:boundsum}
  \expect \exp \Bigl(\eta \sum_{n=0}^N \|w_n\|^2 - \gamma N\Bigr) \le
  \exp \bigl(a\eta \|w_0\|^2\bigr)\;,
\end{equation}    \begin{itemize}\item[]
holds for every $N>0${\rm ,} every $\eta \le \eta_1${\rm ,} and every initial
condition $w_0 \in \CH$.
\ritem{3.}  For every $\eta > 0${\rm ,} there exists a constant
  $C=C(\mathcal{E}_0,\nu,\eta)>0$ such that
  the Jacobian $J_{0,t}$ satisfies almost surely \end{itemize}\vglue-20pt
\begin{equation} \label{e:boundJ}
\|\J_{0,t}^{}\| \le \exp\Bigl(\eta \int_0^t \|w_s\|_1^2\,ds +
Ct\Bigr)\;,
\end{equation}    \begin{itemize}\item[]
for every $t > 0$.
\ritem{4.}  For every $\eta>0$ and every $p > 0${\rm ,} there
  exists $C=C(\mathcal{E}_0,\nu,\eta,p)>0$ such that the Hessian satisfies  
$$
\expect \|\K_{s,t}\|^{p} \le C \exp \bigl(\eta \|w_0\|^2\bigr)\;,
$$
for every $s > 0$ and every $t \in (s,s+1)$.
\end{itemize}
\end{lemma}
 
The proof of \lem{lem:mainbounds} is postponed to
Appendix~A.  \vglue6pt

We now show how to modify the discussion in
Section~\ref{sec:wishfulThinking} to make use of the pathwise
contractivity on small scales to remove the need for the Malliavin
covariance matrix to be invertible on all of $\tilde \CH$.

The point is that since the Malliavin matrix is not invertible, we are
not able to construct a $v\in \L^2([0,T],\R^m)$ for a fixed value of
$T$ that produces the same infinitesimal shift in the solution as an
(arbitrary but fixed) perturbation $\xi$ in the initial condition.
Instead, we will construct a $v\in \L^2([0,\infty),\R^m)$ such that an
infinitesimal shift of the noise in the direction $v$ produces
\textit{asymptotically} the same effect as an infinitesimal
perturbation in the direction $\xi$. In other words, one has
$\|J_{0,t} \xi - A_{0,t} v_{0,t}\| \to 0$ as $t \to \infty$, where
$v_{0,t}$ denotes the restriction of $v$ to the interval $[0,t]$.

Set $\err_t=J_{0,t} \xi - A_{0,t} v_{0,t}$, the residual error for the
infinitesimal variation in the Wiener path $W$ given by $v$. Then we
have from \eref{e:integrationByParts} the \textit{approximate}
integration by parts formula:
\begin{eqnarray}
\scal{\nabla\CP_{t} \phi(w), \xi} &=&\expect_w
\Bigl(\scal{\nabla\bigl(\phi(w_{t})\bigr),\xi}\Bigr)=\expect_w
\Bigl(\bigl(\nabla\phi\bigr)(w_{t})\J_{0,t}\xi\Bigr)\label{e:boundDphi}\\
&=&\expect_w \Bigl(\bigl(\nabla\phi\bigr)(w_{t}) A_{0,t} v_{0,t}\Bigr)
+ \expect_w \bigl(\bigl(\nabla\phi\bigr)(w_{t})\err_{t}\bigr)\nonumber\\
&=&\expect_w\Bigl(\cD^{v_{0,t}} \phi(w_{t})\Bigr) + \expect_w
\bigl(\bigl(\nabla\phi\bigr)(w_{t})\err_{t}\bigr)\nonumber\\ &=&\expect_w
\Bigl(\phi(w_{t}) \int_0^{t} v(s)\,dW(s)\Bigr) + \expect_w
\bigl(\bigl(\nabla\phi\bigr)(w_{t})\err_{t}\bigr)\nonumber\\ &\le&
\|\phi\|_\infty \expect_w \Bigl|\int_0^{t} v(s)\,dW(s)\Bigr| +
\|\nabla \phi\|_\infty \expect_w \|\err_{t}\| \;.\nonumber
\end{eqnarray}
This formula should be compared with \eqref{eq:intbyParts}.
Again if the process $v$ is not adapted to the filtration generated by the increments
of the Wiener process
$W(s)$, the integral must be taken to be a Skorokhod integral;
otherwise It\^o integration can be used.
Note that the residual error satisfies the equation
\begin{equation} \label{e:equrho}
\d_t \err_t = \nu \Delta \err_t + \tilde B(w_t, \err_t) - Qv(t)\;,\quad \err_0 = \xi\;,
\end{equation}   
which can be interpreted as a control problem, where $v$ is the control and 
$\|\err_t\|$ is the quantity that one wants to drive to $0$.

If we can find a $v$ so that $\err_t \rightarrow 0$ as $t \to \infty$ and $ \expect
\bigl|\int_0^{\infty} v(s)\,dW(s)\bigr| < \infty$ then
\eqref{e:boundDphi} and Proposition \ref{l:asfLipFunctions} would imply that $w_t$
is \qsf/. A natural way to accomplish this would be to take
$v(t)=Q^{-1} \tilde B(w_t, \err_t)$, so that $\d_t \err_t = \nu \Delta \err_t$
and hence $\err_t \rightarrow 0$ as $t \rightarrow \infty$. However for
this to make sense it would require that $\tilde B(w_t, \err_t)$ takes
values in the range of $Q$. If the number of Brownian motions $m$ is
finite this is impossible. Even if $m=\infty$, this is still a
delicate requirement which severely limits the range of applicability
of the results obtained (see
\cite{FM}, \cite{b:Fe97}, \cite{b:MattinglySuidan05}).

To overcome these
difficulties, one needs to better incorporate the pathwise smoothing
which the dynamics possesses at small scales. 
Though our ultimate goal is to prove Theorem \ref{theo:main}, which 
 covers \eqref{e:vort} in a fundamentally hypoelliptic
setting, we begin with what might be called the ``essentially
elliptic'' setting. This allows us to outline the ideas in a simpler setting.

\Subsec{Essentially elliptic setting}
\label{sec:EssEllip}
To help to clarify the techniques used in the sections which follow and to
demonstrate their applications, we sketch the proof of the following
proposition which captures the main results of the earlier works on ergodicity,
translated into the framework of the present paper.

\begin{proposition}\label{effectivelyEllipticASF}
Let $\CP_t$ denote the semigroup generated by the solutions to \eref{e:vort} on $\CH$. 
There exists an $N_* = N_*(\CE_0,\nu)$ such that if 
$\CZ_0$ contains $\{k \in \Z^2\;,\; 0<|k| \leq N_*\}${\rm ,} then for any $\eta>0$ there
  exist positive constants $c$ and $\gamma$ so that
  \begin{align*}
       |\nabla \CP_{t}\phi(w)| \leq c\exp\big(\eta \|w\|^2\big)\Big(
        \|\phi\|_\infty + e^{-\gamma t} \|\nabla \phi\|_\infty \Big)\;.
  \end{align*}
\end{proposition}

This result translates the ideas in
\cite{b:EMattinglySinai00}, \cite{MatNS}, \cite{HExp02} to our present
setting. (See also \cite{b:Mattingly03Pre} for more discussion.) The result
does differ from the previous analysis in that it proceeds
infinitesimally. However, both approaches lead to proving the system has a unique
ergodic invariant measure.

The condition on the range of $Q$ can be understood as a type of
``effective ellipticity.'' We will see that the dynamics is
contractive for directions orthogonal to the range of $Q$. Hence
if the noise smooths in these directions, the dynamics will smooth in
the other directions. What directions are contracting depends
fundamentally on a scale set by the balance between $\mathcal{E}_0$
and $\nu$ (see \cite{b:EMattinglySinai00,b:Mattingly03Pre}). 
\prop{l:theWholeEnchilada} holds given a minimal
nondegeneracy condition independent of the viscosity $\nu$, while
\prop{effectivelyEllipticASF} requires a nondegeneracy condition which
depends on $\nu$. 

\demo{Proof of Proposition {\rm \ref{effectivelyEllipticASF}}}
Let $\pi_h$ be the orthogonal projection onto the span of $\{ f_k
: |k| \geq N \}$ and $\pi_\ell=1-\pi_h$. We will fix $N$
presently; however, we will proceed assuming
$\mathcal{H}_{\ell}\eqdef\pi_\ell \mathcal{H} \subset \mbox{Range}(Q)$ and
that $Q_\ell\eqdef\pi_\ell Q$ is invertible on $\mathcal{H}_{\ell}$. By \eref{e:equrho} we therefore 
have full 
control on the evolution of $\pi_\ell \err_t$ by choosing $v$ appropriately.
This allows for an ``adapted'' approach which does not require the control $v$ to
use information about the future increments of the noise process $W$. 

Our approach is first to define a process $\zeta_t$ with the property
that $\pi_\ell \zeta_t$ is $0$ after a finite time and $\pi_h \zeta_t$
evolves according to the linearized evolution, and then choose $v$
such that $\err_t = \zeta_t$. Since $\pi_\ell \zeta_t = 0$ after some
time and the linearized evolution contracts the high modes
exponentially, we readily obtain the required bounds on moments of
$\err_t$.  One can in fact pick any dynamics which are convenient for
the modes which are directly forced. In the case when all of the modes
are forced, the choice  $\zeta_t = (1-t/T)J_{0,t}\xi$ for
$t \in [0,T]$ produces the well-known Bismut-Elworthy-Li formula
\cite{Xuemei}. However, this formula cannot be applied in the present
setting as all of the modes are not necessarily forced.

For $\xi \in
\mathcal{H}$ with $\|\xi\|=1$, define
$\zeta_t$ by
\begin{equation}
  \d_t \zeta_t= -\frac12\frac{\pi_\ell\zeta_t}{\|\pi_\ell\zeta_t\|} 
  +\nu \Delta \pi_h
  \zeta_t + \pi_h \tilde B(w_t,  \zeta_t)\;,\quad \zeta_0=\xi \;.
\end{equation}
With the convention that $0/0 = 0$,
set $\zeta_t^h=\pi_h\zeta_t$ and $\zeta_t^\ell=\pi_\ell\zeta_t$. We
define the infinitesimal perturbation $v$ by 
\begin{equation} \label{e:defvElliptic}
  v(t) =Q_\ell^{-1} F_t\;,\quad
  F_t = \frac12\frac{\zeta_t^\ell}{\|\zeta_t^\ell\|} +
  \nu \Delta \zeta_t^\ell + \pi_\ell \tilde B(w_t,\zeta_t)\;.
\end{equation}   
Because $F_t \in \mathcal{H}_{\ell}$, $Q_\ell^{-1} F_t$ is well defined.  
It is clear from \eref{e:equrho} and \eref{e:defvElliptic} that $\err_t$ and $\zeta_t$
satisfy the same equation, so that indeed $\err_t = \zeta_t$.
Since $\zeta_t^\ell$ satisfies  $\d_t \|\zeta_t^\ell\|^2 = -\|\zeta_t^\ell\|$,
one has $\|\zeta_t^\ell\| \leq \|\zeta_0^\ell\| \leq \|\xi\|=1$. Furthermore, for any initial 
condition $w_0$ and any $\xi$ with $\|\xi\| = 1$, one has $\|\zeta_t^\ell\|=0$ for 
$t \geq 2$.  By calculations similar to those in Appendix~A, there 
exists a constant $C$ so that for any $\eta >0$
\begin{align*}
  \d_t \|\zeta_t^h\|^2 \leq -\Big(\nu N^2 - \frac{C}{\nu \eta^2} -
  \eta\|w_t\|^2_1\Big) \|\zeta_t^h\|^2 +
  \frac{C}{\nu}\|w_t\|_1^2\|\zeta_t^\ell\|^2 \;.
\end{align*}
Hence,
\begin{eqnarray*}
  \|\zeta_t^h\|^2 &\leq& \|\zeta_0^h\|^2 \exp\left( -\Bigl[\nu N^2-
    \frac{C}{\nu \eta^2}\Bigr]t+ \eta\int_0^t \|w_s\|^2_1 ds\right) \\
  &&+C \exp\left( -\Bigl[\nu N^2- \frac{C}{\nu
      \eta^2}\Bigr]\Bigl[t-2\Bigr]+ \eta\int_0^t \|w_r\|^2_1 dr\right)
  \int_0^2\|w_s\|_1^2 ds\;.
\end{eqnarray*}
By \lem{lem:mainbounds}, for any $\eta>0$ and $p\geq 1$ there exist positive constants
$C$ and $\gamma$ so that for all $N$ sufficiently large
\begin{equation}\label{eq:effectivelyElliptic}
  \EE \|\zeta_t^h\|^{p} \leq C(1+\|\zeta_0^h\|^{p}) e^{\eta
  \|w_0\|^2} e^{ -\gamma t} = 2Ce^{\eta
  \|w_0\|^2} e^{ -\gamma t}\;.
\end{equation}
It remains to get control over the size of the perturbation $v$. Since
$v$ is adapted to the Wiener path,
\begin{align*}
  \left(\expect \Bigl|\int_0^{t} v(s)\,dW(s)\Bigr|\right)^2 &\leq
  \int_0^{t}\EE\|v(s)\|^2\,ds\leq C \int_0^t \EE\|F_s\|^2\,ds\;. 
\end{align*}
Now since $\|\pi_\ell\tilde B(u,w)\| \leq C\|u\|\|w\|$ (see
\cite[Lemma A.4]{b:EMattinglySinai00}),
$\|\zeta_t^\ell\|\leq 1$ and $\|\zeta_t^\ell\|=0$ for $t \geq
2$, we see from \eref{e:defvElliptic} that there exists a $C=C(N)$
such that for all $s\geq 0$
$$
   \EE\|F_s\|^2 \leq C \Bigl(1_{\{s \le 2\}} + \EE \|w_s\|^4 \EE \|\zeta_s\|^4 \Bigr)^{1/2}\;.
$$
By using \eqref{eq:effectivelyElliptic} with $p=4$, Lemma
\ref{lem:expbound} from the appendix, to control $\EE \|w_s\|^4$, and picking $N$
sufficiently large, we obtain that for any $\eta>0$ there is a
constant $C$ such that 
\begin{align}
  \label{eq:effEllipFiniteControl}
  \expect \Bigl|\int_0^{\infty} v(s)\,dW(s)\Bigr| \leq
  C \exp\bigl(\eta \|w_0\|^2\bigr)\;.
\end{align}
Since \eqref{eq:effectivelyElliptic} and
\eqref{eq:effEllipFiniteControl} plug into \eqref{e:boundDphi}, the result follows.
\Endproof\vskip4pt  
\Subsec{The truly hypoelliptic setting\/{\rm :}\/ Proof of Proposition~{\rm \ref{l:theWholeEnchilada}}}
\label{sec:proof} 
We now turn to the truly hypoelliptic setting. Unlike  the previous
section, we allow for unstable directions which are not
directly forced by the noise. However, \prop{prop:gen} shows that the 
randomness can reach all of the unstable modes of
interest, i.e.\  those in $\tilde\CH$. In order to show \eref{e:mainbound},
we fix from now on $\xi \in \tilde\CH$ with $\|\xi\| = 1$ and we obtain bounds on
$\scal{\nabla \CP_n \phi(w),\xi}$ that are independent of~$\xi$.

The basic structure of the argument is the same as in the preceding
section on the essentially elliptic setting. We will construct an
infinitesimal perturbation of the Wiener path over the time interval
$[0,t]$ to approximately match the effect on the solution $w_t$ of an 
infinitesimal perturbation
of the initial condition in an arbitrary direction $\xi \in \tilde \CH$.

However, since not all of the unstable directions are in the range of
$Q$, we can no longer infinitesimally correct the effect of the
perturbation in the low mode space as we did in
\eqref{e:defvElliptic}. We rather proceed in a way similar to the
start of Section \ref{sec:wishfulThinking}. However, since the
Malliavin matrix is not invertible, we will regularize it and thus
construct a $v$ which compensates for the perturbation $\xi$ only
asymptotically as $t\to\infty$.  Our construction produces a $v$ which
is \textit{not adapted} to the Brownian filtration, which complicates
a little bit the calculations analogous to
\eqref{eq:effEllipFiniteControl}. A more fundamental difficulty is
that the Malliavin matrix is not invertible on any space which is
easily identifiable or manageable, certainly not on $\L^2_0$. Hence,
the way of constructing $v$ is not immediately obvious.

The main idea for the construction of $v$ is to work with a
regularized version $\tM_{s,t} \eqdef M_{s,t} + \beta$ of the
Malliavin matrix $M_{s,t}$, for some very small parameter $\beta$ to
be determined later. The resulting $\tM^{-1}$ will be an inverse ``up
to a scale'' depending on $\beta$. By this we mean that $\tM^{-1}$
should not simply be thought of as an approximation of $M^{-1}$. It is an
approximation with a very particular form. Theorem~\ref{theo:MP} which is
taken from \cite{b:MattinglyPardoux03Pre} shows that the
eigenvectors with small eigenvalues are concentrated in the small
scales with high probability. This means that $\tM^{-1}$ is very close
to $M^{-1}$ on the large scales\break\vskip-12pt\noindent and very close to the identity times
$\beta^{-1}$ on the small scales. Hence $\tM^{-1}$ will be effective in
controlling the large scales but, as we will see, something else will
have to be done for the small scales.
 
To be more precise, define for integer
values of $n$ the following objects:
$$
\Ja_n = J_{n,n+{1\over 2}}\;,\ \Jb_n = J_{n+{1\over
2},n+1}\;,\ \A_n = A_{n,n+{1\over 2}}\;,\ \M_n=\A_n^{}\A_n^{*}\;,\ \tM_n = \beta +\M_n\;.
$$   
We will then work with a perturbation $v$ which is given by $0$ on all
intervals of the type $[n+{1\over 2},n+1]$, and by $v_n \in
\L^2([n,n+{1\over 2}],\R^m)$ on the remaining intervals. 

We define
the infinitesimal variation $v_n$ by
\begin{equation} \label{e:defv}
v_n = \A_n^*\tM_n^{-1} \Ja_n \err_n\;,
\end{equation}   
where we  denote as before by
$\err_n$ the residual of the infinitesimal displacement at time $n$,
due to the perturbation in the initial condition, which has not yet
been compensated by $v$; i.e.\  $\err_n = J_{0,n} \xi - A_{0,n} v_{0,n}$.
From now on, with a slight abuse of notation we will  write $v_n$
for the perturbation of the Wiener path on $[n,n+{1\over 2}]$ and its
extension (by $0$) to the interval $[n,n+1]$. 

We claim that it follows from \eref{e:defv} that $\err_n$ is given recursively by
\begin{equation} \label{e:defxin}
\err_{n+1} = \Jb_n \beta \tM_n^{-1} \Ja_n \err_n\;,
\end{equation}   
with $\err_0 = \xi$.  To see the claim observe that \eref{e:defxin} implies $\J_{n,n+1}
\err_n = \Jb_n \Ja_{n} \err_n = \Jb_n \A_n v_n + \err_{n+1}$.  Using
this and the definitions of the operators involved, we see that indeed
\begin{align*}
  A_{0,N} v_{0,N}&=\sum_{n=0}^{N-1} \J_{(n+1),N} \Jb_n \A_n
  v_n=\sum_{n=0}^{N-1} \bigl(\J_{n,N} \err_n - \J_{(n+1),N}
  \err_{n+1}\bigr) \\&=\J_{0,N}\xi - \err_{N} \; .
\end{align*}
Thus, at time $N$, the infinitesimal variation in the
Wiener path $v_{0,N}$ corresponds to the infinitesimal perturbation in
the initial condition $\xi$ up to an error $\err_{N}$.

It therefore remains to show that this choice of $v$ has desirable
properties.  In particular we need to demonstrate properties similar
to \eqref{eq:effectivelyElliptic} and
\eqref{eq:effEllipFiniteControl}. The analogous statements are given
by the next two propositions whose proofs will be the content of Sections
\ref{sec:hypoError} and \ref{sec:costControl}. Both of these propositions rely
heavily on the following theorem obtained in \cite[Th.~6.2]{b:MattinglyPardoux03Pre}.

\begin{theorem}\label{theo:MP}
  Denote by $M$ the Malliavin matrix over the time interval
  $[0,{1\over 2}]$ and define $\tilde \CH$ as above.  For every
  $\alpha,\eta,p$ and every orthogonal projection $\pi_\ell$ on a
  finite number of Fourier modes{\rm ,} there exists $\tilde C$ such that
\begin{equation} \label{e:lemMP}
  \P \bigl(\scal{M\phi,\phi} < \eps \|\phi\|_1^2\bigr) \le \tilde C
  \eps^p\exp \bigl(\eta \|w_0\|^2\bigr)\;,
\end{equation}   
holds for every \/{\rm (}\/random\/{\rm )}\/ vector $\phi \in \tilde \CH$ satisfying $\|\pi_\ell \phi\|
\ge \alpha\|\phi\|_1$ almost surely{\rm ,} for every $\eps \in (0,1)${\rm ,} and for every $w_0 \in \tilde\CH$.
\end{theorem}

The next proposition shows that we can construct a $v$ which has the desired effect
of driving the error
$\err_t$ to zero as $t \rightarrow \infty$. 

\begin{proposition}\label{prop:boundxi}
For any $\eta >0${\rm ,} there exist constants $\beta > 0$ and $C>0$ such that
\begin{equation} \label{e:bounderr}
\expect \|\err_N\|^{10} \le {C\exp(\eta\|w_0\|^2) \over 2^N} 
\end{equation}   
holds for every $N > 0$. \/{\rm (}\/Note that by increasing $\beta$ further{\rm ,} 
 the $2^N$ in the denominator could be replaced by 
$K^N$ for an arbitrary $K \ge 2$ without altering the value of $C$.{\rm )}
\end{proposition}

However for the above result to be
useful, the ``cost'' of shifting the noise by $v$ (i.e.\  the norm of $v$ in the Cameron-Martin
space) must be finite. Since the time
horizon is infinite, this is not a trivial requirement. In the
``essentially elliptic'' setting, it was demonstrated in
\eqref{eq:effEllipFiniteControl}. In the ``truly hypoelliptic'' setting, we obtain

\begin{proposition}\label{prop:hypoCost} For any $\eta >0${\rm ,} there exists a
  constant $C$ so that  
  \begin{equation}\label{e:hypoCost} 
    \EE \Bigl|\int_0^{N} v_{0,s}\,dW(s)\Bigr|^2 \leq \frac{C}{\beta^2}
    e^{\eta\|w_0\|^2} \sum_{n=0}^\infty \bigl(\expect
    \|\err_n\|^{10}\bigr)^{\frac15} .
  \end{equation}
\/{\rm (}\/Note that the power $10$ in this expression is arbitrary and can be
brought as close to $2$ as one wishes.{\rm )}
\end{proposition}

Plugging these estimates into \eqref{e:boundDphi}, we obtain
Proposition~\ref{l:theWholeEnchilada}.
Note that even though Proposition~\ref{l:theWholeEnchilada} is sufficient for
the present article, small modifications of \eqref{e:boundDphi} produce the following stronger bound.

\begin{proposition}
For every $\eta > 0$ and every $\gamma > 0${\rm ,} there exist 
constants $C_{\eta,\gamma}$ such
that for every Fr\ee chet differentiable function $\phi$ from $\tilde \CH$
to~$\R${\rm ,}
$$
\|\nabla\CP_n \phi(w)\| \le \exp(\eta\|w\|^2)\Bigl(C_{\eta,\gamma} \sqrt{\bigl(\CP_n |\phi|^2\bigr)(w)}
      +  \gamma^n \sqrt{\bigl(\CP_n\|\nabla\phi\|^2\bigr)(w)}\Bigr)\;,
$$ 
for every $w \in \tilde \CH$ and $n \in \N$.
\end{proposition}

\Proof 
Applying Cauchy-Schwarz to the terms on the right-hand side of the penultimate line of 
\eref{e:boundDphi} one obtains
$$
|\scal{\nabla\CP_n \phi,\xi}| \le \Bigl(\EE \Bigl|\int_0^{n} v_{0,s}\,dW(s)\Bigr|^2 \CP_n |\phi|^2\Bigr)^{1/2}
      +  \Bigl(\expect \|\err_n\|^{10} \CP_n\|\nabla\phi\|^2\Bigr)^{1/2}\;.
$$
It now suffices to use the bounds from the above propositions
and to note that the right-hand side is independent of the choice of $\xi$ provided $\|\xi\| = 1$.
\hfq

\Subsec{Controlling the error\/{\rm :}\/ Proof of Proposition {\rm \ref{prop:boundxi}}}
\label{sec:hypoError}
Before proving \prop{prop:boundxi}, we state the following lemma, which summarizes
the effect of our control on the perturbation and will be proved at the end of this section.
\begin{lemma}\label{lem:boundxi2}
  For every two constants $\gamma, \eta > 0$ and every $p\ge 1${\rm ,} there
  exists a constant $\beta_0 > 0$ such that
$$
  \expect\bigl(\|\err_{n+1}\|^p\,|\, \CF_n\bigr) \le \gamma e^{\eta
    \|w_n\|^2} \|\err_n\|^p
$$  
holds almost surely whenever $\beta \le \beta_0$.
\end{lemma}
 
\demo{Proof of \prop{prop:boundxi}}
Define
$$
C_n = {\|\err_{n+1}\|^{10} \over \|\err_n\|^{10}} \;,
$$
with the convention that $C_n = 0$ if $\err_n = 0$. Note that since $\|\err_0\| = 1$,
one has $\|\err_N\|^{10} = \prod_{n=0}^{N-1} C_n$. 
We begin by
establishing some properties of $C_n$ and then use them to prove the proposition.

Note that $\|\beta \tM_n^{-1}\| \leq 1$ and so, by
\eref{e:boundJ} and \eref{e:defxin}, for every $\eta > 0$
there exists a constant $C_\eta > 0$ such that
\begin{equation} \label{e:aprioriCn}
C_n \le \|\Jb_n \beta \tM_n^{-1} \Ja_n\|^{10} \le \|\Jb_n\|^{10} \|\Ja_n\|^{10} \le \exp \Bigl({\eta \int_n^{n+1} \|w_s\|_1^2\,ds + C_\eta}\Bigr)\;,
\end{equation}   
almost surely. Note that this bound is independent of  $\beta$.  
Next, for given values of $\eta$ and $R > 0$, we define
$$
  C_{n,R} = \left\{\begin{array}{ll} {e^{-\eta R}} & \text{if
        $\|w_n\|^2 \ge 2R$,} \\[0.5em] e^{\eta R}C_n & \text{otherwise.}
    \end{array}\right.
$$   
Obviously both $C_n$ and $C_{n,R}$ are $\CF_{n+1}$-measurable.
\lem{lem:boundxi2} shows that for every $R> \eta^{-1}$, one can find a
$\beta > 0$ such that
\begin{equation} \label{e:asbound}
\expect \bigl(C_{n,R}^2\,|\,\CF_n\bigr) \le {1\over 2}\;,\quad
\text{almost surely}\;.
\end{equation}   
Note now that \eref{e:aprioriCn} and the definition of $C_{n,R}$
immediately imply that
\begin{equation} \label{e:boundCn}
  C_n \le C_{n,R} \exp \Bigl({\eta \int_n^{n+1} \|w_s\|_1^2\,ds + \eta
    \|w_n\|^2 + C_\eta - \eta R}\Bigr)\;,
\end{equation}   
almost surely. 
This in turn implies that
\begin{eqnarray*}
  \prod_{n=0}^{N-1} C_n &\le& \prod_{n=0}^{N-1} C_{n,R}^2 + \prod_{n=0}^{N-1} \exp
  \Bigl({2\eta \int_n^{n+1} \|w_s\|_1^2\,ds + 2\eta \|w_n\|^2 + 2C_\eta
    - 2\eta R}\Bigr)\\
  &\le& \prod_{n=0}^{N-1} C_{n,R}^2 + \exp\Bigl({4\eta\sum_{n=0}^{N-1}
    \|w_n\|^2 + 2N(C_\eta - \eta R)}\Bigr)\\
  &&+ \exp \Bigl({4\eta \int_0^{N} \|w_s\|_1^2\,ds + 2N(C_\eta -
    \eta R)}\Bigr)\;.
\end{eqnarray*}
Now fix $\eta > 0$ (not too large). In light of \eqref{e:boundInt} and
\eqref{e:boundsum}, we can then choose $R$ sufficiently large so that
the  last  two terms satisfy the required bounds. Then, we choose
$\beta$ sufficiently small so that \eref{e:asbound} holds and the
estimate follows.
\Endproof\vskip4pt  

To prove Lemma \ref{lem:boundxi2}, we will use the following two
lemmas. The first is simply a consequence of the dissipative nature of the
equation. Because of the Laplacian, the small scale perturbations are
strongly damped.
\begin{lemma}\label{l:HighContract}
  For every $p\ge 1${\rm ,} every $T>0${\rm ,} and every two constants $\gamma,
  \eta > 0${\rm ,} there exists an orthogonal projector $\pi_\ell$ onto a finite number of
  Fourier modes such that
\begin{align}  \expect \|(1-\pi_\ell^{}) \J_{0,T}^{}\|^{p} &\le \gamma
  \exp\bigl(\eta \|w_0\|^2\bigr)\;,\label{e:boundPhJ1}\\[7pt] \expect
  \|\J_{0,T}^{}(1-\pi_\ell^{})\|^{p} &\le \gamma \exp\bigl(\eta
  \|w_0\|^2\bigr)\;.\label{e:boundPhJ2}
\end{align}
\end{lemma}
\vglue8pt

The proof of the above lemma is postponed to the appendix. The
second lemma is central to the hypoelliptic results in this paper.
It is the analog of \eqref{eq:effEllipFiniteControl} from the
essentially elliptic setting and provides the key to controlling the
``low modes'' when they are not directly forced and Girsanov's theorem
cannot be used directly. This result makes use of the results in
\cite{b:MattinglyPardoux03Pre} which contains the heart of the
analysis of the structure of the Malliavin matrix for equation
\eqref{e:vort} in the  hypoelliptic setting.

\begin{lemma}\label{lem:nicebound}
Fix $\xi \in \tilde \CH$ and define
$$
\zeta = \beta (\beta + \M_0)^{-1} \hat J_{0} \xi\;.
$$  
Then{\rm ,} for every two constants $\gamma, \eta > 0$ and every low-mode
orthogonal projector~$\pi_\ell${\rm ,} there exists a constant $\beta > 0$ such that
$$
\expect \|\pi_\ell\zeta\|^p \le \gamma e^{\eta \|w_0\|^2} \|\xi\|^p\;.
$$   
\end{lemma}
\vglue8pt

\begin{remark}\label{rem:choicexi}
  Since one has obviously that $\|\zeta\|Ê\le \|\hat J_0\xi\|$, this lemma
  tells us that applying the operator $\beta (\beta + \M_0)^{-1}$ (with a
  very small value of $\beta$) to a vector in $\tilde\CH$ either reduces its norm
  drastically or transfers most of its ``mass'' into the high modes
  (where the cutoff between ``high'' and ``low'' modes is arbitrary
  but influences the possible choices of $\beta$). This explains why
the control $v$ is set to $0$ for half of the time in Section~\ref{sec:proof}.  In order to 
ensure that the norm of $\err_n$ gets really reduced after one step, we choose the control
in such a way that $\beta (\beta + \M_n)^{-1} \Ja_n$ is composed by $\Jb_n$, using the fact
  embodied in Lemma \ref{l:HighContract} that the Jacobian will
  contract the high modes before the low modes start to grow out of control.
\end{remark}

{\it Proof of Lemma {\rm \ref{lem:nicebound}}}.
For $\alpha > 0$, let $A_\alpha$ denote the event $\|\pi_\ell \zeta\| > \alpha \|\zeta\|_1$.
  We also define the random vectors
$$
\zeta_\alpha(\omega) = \zeta(\omega) \chi_{A_\alpha}(\omega)\;,\quad
\bar \zeta_\alpha(\omega) = \zeta(\omega)-\zeta_\alpha(\omega)\;,\quad \omega \in \Omega\;,
$$ 
where  $\omega$ is the chance variable and $\chi_A$ is the
characteristic function of a set $A$.  It is clear that
$$
  \expect \|\pi_\ell\bar \zeta_\alpha\|^p \le \alpha^p \,\expect
  \|\zeta\|_1^p\;.
$$ 
Using the bounds \eref{e:boundJ} and \eref{e:boundInt} on the Jacobian
and the fact that $\M_0$ is a bounded operator from $\CH_1$ (the Sobolev space
 \pagebreak of functions with square integrable derivatives) into $\CH_1$,
we get
\begin{equation} \label{e:boundddd}
  \expect \|\pi_\ell\bar \zeta_\alpha\|^p \le \alpha^p \expect
  \|\zeta\|_1^p \le \alpha^p \expect \|\Ja_0\xi\|_1^p \le {\gamma\over 2}
  e^{\eta \|w_0\|^2} \|\xi\|^p\;,
\end{equation}   
(with $\eta$ and $\gamma$ as in the statement of the proposition) 
for sufficiently small $\alpha$. From now on, we fix $\alpha$ such that \eref{e:boundddd} holds.
One has the chain of inequalities
\begin{eqnarray} \label{e:chain}
  \scal{\zeta_\alpha, \M_0\zeta_\alpha} &\le& \scal{\zeta, \M_0\zeta} \le
  \scal{\zeta, (\M_0+\beta)\zeta} \\ &=& \beta \scal{\Ja_0\xi,
    \beta(\M_0+\beta)^{-1}\Ja_0\xi} \le \beta \|\Ja_0\xi\|^2\;. \nonumber
\end{eqnarray}
From \theo{theo:MP}, we see  furthermore  that, for every $p_0 > 0$, there exists
a constant $\tilde C$ such that
$$
  \P \bigl(\scal{\M_0 \zeta_\alpha, \zeta_\alpha} < \eps
  \|\zeta_\alpha\|_1^2\bigr) \le \tilde C \eps^{p_0}\exp \bigl(\eta
  \|w_0\|^2\bigr)\;,
$$
holds for every $w_0 \in \tilde \CH$ and every $\eps \in (0,1)$.
Consequently, 
$$
  \P \Bigl({\|\zeta_\alpha\|_1^2\over \|\Ja_0\xi\|^2} > {1 \over \eps}
  \Bigr) \le 
 \P \bigl(\scal{\M_0 \zeta_\alpha, \zeta_\alpha} < \eps\beta
  \|\zeta_\alpha\|_1^2\bigr) \le
  \tilde C \beta^{p_0} \eps^{p_0}\exp \bigl(\eta \|w_0\|^2\bigr)\;,
$$
where we made use of \eref{e:chain} to get the first inequality.
This implies that, for every $p,q \ge 1$, there exists a constant $\tilde C$ such that
\begin{equation} \label{e:blablabla}
  \expect \Bigl({\|\zeta_\alpha\|_1^p\over \|\Ja_0\xi\|^p}\Bigr) \le
  \tilde C \beta^q \exp \bigl(\eta \|w_0\|^2\bigr)\;.
\end{equation}   
Since $\|\pi_\ell\zeta_\alpha\| \le \|\zeta_\alpha\|_1$ and
$$
  \expect \|\zeta_\alpha\|_1^p \le \sqrt{\expect
    \Bigl({\|\zeta_\alpha\|_1^{2p}\over \|\Ja_0\xi\|^{2p}}\Bigr) \expect
    \|\Ja_0\xi\|^{2p}}\;,
$$ 
it follows from \eref{e:blablabla} and the bound \eref{e:boundJ} on
the Jacobian that, by choosing $\beta$ sufficiently small, one gets
\begin{equation} \label{e:boundd}
  \expect \|\pi_\ell \zeta_\alpha\|^p \le {\gamma\over 2} e^{\eta
    \|w_0\|^2} \|\xi\|^p\;.
\end{equation}   
Note that $\expect \|\pi_\ell\zeta\|^p
= \expect \|\pi_\ell \zeta_\alpha\|^p + \expect \|\pi_\ell\bar
\zeta_\alpha\|^p$ since only one of the previous two terms is nonzero 
for any given realization $\omega$. The claim thus follows from \eref{e:boundddd}
and \eref{e:boundd}.
\Endproof\vskip4pt  

Using \lem{l:HighContract} and \lem{lem:nicebound}, we now give the

\demo{Proof  of Lemma {\rm \ref{lem:boundxi2}}}
  Define $\zeta_n = \beta\tM_n^{-1}\Ja_n \err_n$, so that $\err_{n+1} =
  \Jb_n \zeta_n$.  It follows from the definition of $\tM_n$ and
  the bounds \eref{e:boundJ} and \eref{e:boundInt} on the Jacobian
  that there exists a constant $C$ such that
$$
  \expect\bigl(\|\zeta_n\|^p\,|\, \CF_n\bigr) \le C e^{{\eta\over 2}
    \|w_n\|^2} \|\err_n\|^p\;,
$$   
uniformly in $\beta > 0$. Applying \eref{e:boundPhJ2} to this bound
yields the existence of a projector $\pi_\ell$ on a finite number of
Fourier modes such that
$$
  \expect\bigl(\|\Jb_n (1-\pi_\ell) \zeta_n\|^p\,|\, \CF_n\bigr)
  \le \gamma e^{\eta \|w_n\|^2} \|\err_n\|^p\;.
$$ 
Furthermore, \lem{lem:nicebound} shows that, for an arbitrarily small value $\tilde\gamma$,
one can choose $\beta$
sufficiently small so that
$$
  \expect\bigl(\|\pi_\ell \zeta_n\|^p\,|\, \CF_n\bigr) \le \tilde
  \gamma e^{{\eta\over 2} \|w_n\|^2} \|\err_n\|^p\;.
$$
Applying again the \textit{a priori} estimates
\eref{e:boundJ} and \eref{e:boundInt} on the Jacobian, we see that one can choose
$\tilde\gamma$ (and thus $\beta$) sufficiently small so that
$$
  \expect\bigl(\|\Jb_n \pi_\ell \zeta_n\|^p\,|\, \CF_n\bigr) \le
  \gamma e^{\eta \|w_n\|^2} \|\err_n\|^p\;,
$$
and the result follows.
\hfq

\Subsec{Cost of the control\/{\rm :}\/ Proof of Proposition~{\rm \ref{prop:hypoCost}}}
\label{sec:costControl}
Since the process $v_{0,s}$ is not adapted to the Wiener process
$W(s)$, the integral must be taken to be a Skorokhod integral.  We
denote by $\cD_s F$ the Malliavin derivative of a random variable $F$
at time $s$ (see \cite{MR96k:60130} for definitions). Suppressing the
dependence on the initial condition $w$, we obtain from the definition
of the Skorokhod integral and from the corresponding It\^o isometry
(see e.g.\   \cite[p.~39]{MR96k:60130})
$$
{\expect
      \Bigl|\int_0^{N} v(s)\,dW(s)\Bigr|^2}  \le {\expect
      \|v_{0,N}\|^2 + \sum_{n=0}^N \int_{n}^{n+{1\over 2}} \int_n^{n+{1\over 2}} \expect
      \norm{\cD_s v_n(t)}^2\,ds\,dt}\;.
$$
(Remember that $v_n(t) = 0$ on $[n+{1\over 2},n+1]$.)
In this expression, the norm $\norm{\cdot}$ denotes the
Hilbert-Schmidt norm on $m\times m$ matrices, so that  one has
$$
  \int_{n}^{n+{1\over 2}} \int_n^{n+{1\over 2}} \expect \norm{\cD_s v_n(t)}^2\,ds\,dt =
  \sum_{i=1}^m \int_n^{n+{1\over 2}} \expect \|{\cD_s^i v_n^{}}\|^2\,ds\;,
$$   
where the norm $\|{\cdot}\|$ is in $\L^2([n,n+{1\over 2}],\R^m)$ and $\cD_s^i$
denotes the Malliavin derivative with respect to the $i^{\rm th}$ component
of the noise at time $s$.

In order to obtain an explicit expression for $\cD_s^i v_n^{}$, we start by computing separately
the Malliavin derivatives of the various expressions that enter into its construction.
Recall from \cite{MR96k:60130} that $\cD_s^i w_t = J_{s,t} Q e_i$ for $s < t$. It follows
from this and the expression \eref{eq:JLinear} for the Jacobian that the Malliavin derivative
of $J_{s,t}\xi$ is given by
$$
\d_t \cD_r^i J_{s,t}\xi = \nu \Delta \cD_r^i J_{s,t}\xi  + \tilde B(w_t, \cD_r^i J_{s,t}\xi )
+ \tilde B(J_{r,t} Q e_i, J_{s,t} \xi)\;.
$$  
From the variation of constants formula and the expression \eref{e:defK} for the process $K$,
we get
\begin{equation} \label{e:derJac}
\cD_r^i J_{s,t}\xi = \left\{\begin{array}{rl} K_{r,t}(Qe_i, J_{s,r} \xi) & \text{if $r \ge s$,} \\[0.5em]
        K_{s,t}(J_{r,s}Qe_i,\xi)& \text{if $r \le s$.} \end{array}\right.
\end{equation}   
In the remainder of this section, we will use the convention that if $A:\CH_1 \to \CH_2$ is
a random linear map between two Hilbert spaces, we denote by $\cD_s^i A : \CH_1 \to \CH_2$
the random linear map defined by
$$
(\cD_s^i A) h = \scal{\cD_s (A h), e_i}\;.
$$   
With this convention, \eref{e:derJac} immediately yields 
\begin{equation} \label{e:derJn}
\cD_r^i \hat \J_n^{} w = \K_{r,n+{1\over 2}}\bigl(\J_{n,r} w,Qe_i\bigr)\quad \text{for $r \in [n,n+{1\over 2}]$.}
\end{equation}   
Similarly, we see from \eref{e:derJac} and the definition of
$\A_n$ that the map $\cD_r^i \A_n^{}$ given by
\begin{eqnarray}
  \cD_r^i \A_n^{} h &=&   \int_n^r \K_{r,n+{1\over 2}} \bigl(\J_{s,r}Q h(s),Qe_i)\bigr)\,ds
\label{e:DAn}\\ & & + \int_r^{n+{1\over 2}} \K_{r,n+{1\over 2}} \bigl(Q h(s),
\J_{r,s}Qe_i)\bigr)\,ds\;. \nonumber
\end{eqnarray}
We denote its adjoint by $\cD_r^i \A_n^{*}$. Since $\tM_n = \beta + \A_n  \A_n^*$,
we get from the chain rule
$$
\cD_s^i \tM_n^{-1} = - \tM_n^{-1} \Bigl(\bigl(\cD_s^i \A_n^{}\bigr)\A_n^* +
  \A_n^{}\bigl(\cD_s^i \A_n^*\bigr)\Bigr) \tM_n^{-1}\;.
$$ 
Since $\err_n$ is $\CF_n$-measurable, one has $\cD_r^i \err_n = 0$ for
$r \ge n$. Therefore, combining the above expressions with the Leibniz rule
applied to the definition \eref{e:defv} of $v_n$ yields
\begin{eqnarray*}
  \cD_s^i v_n^{} &=& \bigl(\cD_s^i \A_n^*\bigr)\tM_n^{-1} \hat\J_n \err_n +
  \A_n^*\tM_n^{-1} \bigl(\cD_s^i\hat \J_n^{}\bigr) \err_n \\ && -
  \A_n^*\tM_n^{-1}\Bigl(\bigl(\cD_s^i \A_n^{}\bigr)\A_n^* +
  \A_n^{}\bigl(\cD_s^i \A_n^*\bigr)\Bigr)\tM_n^{-1} \hat\J_n \err_n\;.
\end{eqnarray*}
Since $\tM_n = \beta + \A_n^{} \A_n^*$, one has the almost sure bounds
$$
\|\A_n^* \tM_n^{-1/2}\| \le 1\;,\quad \|\tM_n^{-1/2}\A_n\| \le 1\;,\quad \|\tM_n^{-1/2}\|\le \beta^{-1/2}\;.
$$   
This immediately yields
$$
\|{\cD_s^i v_n^{}}\| \le {3\beta^{-1}} \|{\cD_s^i \A_n^{}}\| \|\hat\J_n\|
\|\err_n\| + \beta^{-1/2} \|\cD_s^i \hat\J_n^{}\| \|\err_n\|\;.
$$ 
Combining this with \eref{e:DAn}, \eref{e:derJn}, and \lem{lem:mainbounds}, we obtain,
for every $\eta > 0$, the existence of a constant $C$ such that
$$
\expect \|{\cD_s^i v_n^{}}\|^2 \le C e^{\eta\|w\|^2} \beta^{-2}
\bigl(\expect \|\err_n\|^{10}\bigr)^{1\over 5}\;.
$$
Applying \lem{lem:mainbounds} to the definition of $v_n$ we easily
get a similar bound for $\expect \|v_n^{}\|^2$, which then implies the quoted result.

\section{Discussion and conclusion} 
\label{sec:conclusion}


Even though the results obtained in this work are relatively complete, they still leave a few 
questions open.

Do the transition probabilities for \eref{e:vort} converge towards the invariant measure
and at what rate? In other words, do the solutions to \eref{e:vort} have the mixing property?
We expect this to be the case and plan to answer this question in a subsequent publication.

What happens if $\CH \neq \tilde \CH$ and one starts the system with an initial condition
$w_0 \in \CH \setminus \tilde \CH$? If the viscosity is sufficiently large, we know that the component
of $w_t$ orthogonal to $\tilde \CH$ will decrease exponentially with time. This is however not expected
to be the case when $\nu$ is small. In this case, we expect to have (at least) one invariant measure
associated to every (closed) subspace $V$ invariant under the flow.

\appendix

\setcounter{section}{1}

\vglue16pt \centerline{\bf Appendix A. {\it A priori} estimates for the Navier-Stokes equations} 

\demo{Note} The letter $C$ denotes generic constants whose value can change from one
line to the next even within the same equation. The possible dependence of $C$ on
the parameters of \eref{e:vort} should be clear from the context.\Enddemo

We define for $\alpha \in \R$ and for $w$ a smooth function on
$[0,2\pi]^2$ with mean $0$ the norm $\|w\|_\alpha$ by
$$
\|w\|_\alpha^2 = \sum_{k \in \Z^2\setminus\{0,0\}} |k|^{2\alpha}
w_k^2\;,$$
   where of course $w_k$ denotes the Fourier mode with
wavenumber $k$.  Define furthermore $(\KK w)_k = -i w_k
k^{\perp}/\|k\|^2$, $B(u,v) = (u\cdot\nabla)v$ and 
$\mathcal{S}=\{
s=(s_1,s_2,s_3) \in \R_+^3 : \sum s_i \geq 1, s\neq
(1,0,0),(0,1,0),(0,0,1) \}$. Then the following relations are useful
(cf. \cite{b:CoFo88}):  
\begin{alignat}{5} \scal{B(u,v),w} &= -
\scal{B(u,w),v} \quad && \text{if $\nabla \cdot u = 0$}, \label{e:idB}&\\
|\scal{B(u,v),w}| &\le C\|u\|_{s_1} \|v\|_{1+s_2} \|w\|_{s_3},\quad&&
(s_1,s_2,s_3) \in \mathcal{S}, \label{e:boundB}&\\ \|\KK u\|_\alpha &= 
\|u\|_{\alpha - 1},\label{e:idK} \\ \|u\|_\beta^2 &\le  \eps
\|u\|_\alpha^2 + \eps^{-2{\gamma - \beta \over
    \beta-\alpha}}\|u\|_\gamma^2 \quad && \text{if $0 \le \alpha <
  \beta < \gamma$ and $\eps > 0$.}& \label{e:interp}
\end{alignat}

Before we turn to the proof of \lem{lem:mainbounds}, we give the
following essential bound on the solutions of \eref{e:vort}.

\begin{lemma}\label{lem:expbound}
There exist constants $\eta_0>0$ and $C>0${\rm ,} such that for
every $t>0$ and every $\eta \in (0,\eta_0]${\rm ,} the bound
\begin{equation} \label{e:expbound}
\expect \exp(\eta \|w_t\|^2) \le C \exp(\eta e^{-\nu t}
\|w_0\|^2)
\end{equation}   
holds.
\end{lemma}

\Proof 
From \eref{e:idB} and It\^os formula, we obtain
$$
\|w_t\|^2 - \|w_0\|^2 + 2\nu \int_0^t \|w_r\|_1^2\,dr = \int_0^t
\scal{w_r, Q\,dW(r)} + \mathcal{E}_0 t\;,
$$ 
where  $\mathcal{E}_0 = \tr QQ^*$.  Using the fact that $\|w_r\|_1^2 \ge
\|w_r\|^2$, we get
$$
  \|w_t\|^2 \le e^{-\nu t}\|w_0\|^2 + {\mathcal{E}_0 \over \nu} + \int_0^t
  e^{-\nu(t-r)} \scal{w_r, Q\,dW(r)} - \nu \int_0^t
  e^{-\nu(t-r)} \|w_r\|^2\,dr\;.
$$
There exists a constant $\alpha > 0$ such that $\nu
\|w_r\|^2 > {\alpha \over 2}\|Q^*w_r\|^2$, so that \cite[Lemma~A.1]{b:Mattingly02d}
implies
$$
  \P \Bigl(\|w_t\|^2 - e^{-\nu t}\|w_0\|^2 - {\mathcal{E}_0 \over \nu} > {K
    \over \alpha}\Bigr) \le e^{-K}\;.
$$
Note now that if a random variable $X$ satisfies $\P(X \ge C) \le 1/C^2$ for
all $C \ge 0$, then $\expect X \le 2$.
The bound \eref{e:expbound} thus follows immediately with for example $\eta_0 = \alpha/2$ and
$C = 2 \exp({\alpha \CE_0\over 2\nu})$.
\Endproof\vskip4pt  

We now turn to the proof of  \lem{lem:mainbounds}.
\demo{Proof  of \lem{lem:mainbounds}}
   {\it Point} 1. From \eref{e:idB} and It\^o's formula, for any $\eta > 0$ we obtain
\begin{multline*}
\eta \|w_t\|^2+ \eta \nu\int_s^t\|w_r\|_1^2 \,dr -\eta \mathcal{E}_0 (t-s)\;\\=  \eta \|w_s\|^2
+ \eta \int_s^t \scal{w_r, Q\,dW(r)} - \eta \nu\int_s^t \|w_r\|_1^2 \,dr 
\end{multline*}
where  $\mathcal{E}_0 = \tr QQ^*$.  Denote by $M(s,t)$ the first two
terms on the right-hand side of the last expression and set
$N(s,t)=M(s,t)- \eta \nu\int_s^t\|w_r\|_1^2 \,dr$. Now observe that
with  $\alpha$ as in the proof of \lem{lem:expbound} above,
one has $N(s,t) \leq M(s,t) -
\frac{\alpha}{2\eta} \langle M \rangle(s,t)$ where $\langle M
\rangle(s,t)$ is the quadratic variation of the continuous $L^2$-martingale $M$. Hence
by the standard exponential martingale estimate, $\P( \sup_{t\geq s}
N(s,t) \geq K\,|\,\CF_s) \leq \exp\,(\eta \|w_s\|^2-{\alpha K\over \eta})$ for all $s\geq0$.  
Here we use the notation $\mathcal{F}_s$ to denote the filtration generated by the noise up to
the time $s$. Thus, for all $\eta \in (0,\alpha/2]$ and $s\geq 0$,
$$
  \expect \exp\Bigl( \eta \sup_{t \geq s}\Bigl(\|w_t\|^2+ \nu  \int_s^t
  \|w_r\|_1^2\,dr -\mathcal{E}_0 (t-s)\Bigr)\,\Big|\, \mathcal{F}_s\Bigr) 
\leq 2  \exp \bigl(\eta \|w_s\|^2\bigr)\;.
$$  
Choosing $\eta_0$ as above and using  Lemma \ref{lem:expbound} to bound the expected
value of the right-hand side, we complete the proof.
 
\demo{Point {\rm 2}} Taking conditional expectations with
respect to $\CF_{N-1}$ on the left hand side of \eref{e:boundsum} and
applying \lem{lem:expbound}, one has
$$
  \expect \exp \Bigl(\eta \sum_{n=0}^N \|w_n\|^2\Bigr) \le C\expect
  \exp \Bigl(\eta e^{-\nu} \|w_{N-1}\|^2 + \eta \sum_{n=0}^{N-1}
  \|w_n\|^2\Bigr)\;.
$$ 
Applying this procedure repeatedly, one obtains
$$
  \expect \exp \Bigl(\eta \sum_{n=0}^N \|w_n\|^2\Bigr) \le C^N \exp(a
  \eta \|w_0\|^2)\;,
$$   
where $a = \sum_{n=0}^\infty e^{-\nu n}$. This computation is
valid, provided $a \eta$ is smaller than $\eta_0$, so the result
follows by taking $\eta_1 = \eta_0 / a$.

\demo{Point {\rm 3}}
We define $\uu_t = \J_{0,t} \uu_0$ for some $\uu_0 \in \CH$. The
evolution of $\uu_t$ is then given by \eqref{eq:JLinear}.
We thus have for the $\CH$-norm of $\xi$ the equation
$$
  \d_t \|\uu_t\|^2 = -2\nu \|\nabla \uu_t\|^2 + 2 \scal{B(\KK
    \uu_t,w_t), \uu_t}\;.
$$  
Equation \eref{e:idB} yields the existence of a constant $C$ such that
$2|\scal{B(\KK h,w), \zeta}| \le C\|w\|_1 \|h\| \|\zeta\|_{1/2}$
for example. By interpolation, we get
\begin{equation} \label{e:boundinterp}
  2|\scal{B(\KK h,w), \zeta}| \le \nu \|\zeta\|_1^2 + {C \over
    \eta^2 \nu} \|\zeta\|^2 + {\eta\over 2} \|w\|_1^2 \|h\|^2\;,
\end{equation}   
and therefore
\begin{equation} \label{e:boundxi1}
  \d_t \|\uu_t\|^2 \le -\nu \|\nabla \uu_t\|^2 + {C \over \eta^2 \nu}
  \|\uu_t\|^2 + {\eta\over 2} \|w_t\|_1^2 \|\uu_t\|^2\;,
\end{equation}   
for every $\eta > 0$. This yields \eref{e:boundJ}.
 
\demo{Point {\rm 4}}
This bound follows in a rather straightforward way from
\eref{e:boundzeta} which is in the next proof. Standard Sobolev estimates and interpolation
inequalities give for the symmetrized nonlinearity $\tilde B$ the
bound
\begin{eqnarray*}
\|\tilde B(u,w)\| &\le &C \bigl(\|u\|_{1/2}\|w\|_1 +
\|u\|_{1}\|w\|_{1/2}\bigr) \\ &\le &C \bigl(\|u\|^{1/2} \|u\|_1^{1/2}
\|w\|_1 + \|w\|^{1/2} \|w\|_1^{1/2} \|u\|_1\bigr) \;.
\end{eqnarray*}
Combining this with the definition \eref{e:defK} of $\K_{s,t}$ and
bound \eref{e:boundzeta} yields for $s,t \in [0,1]$
\begin{eqnarray*}
\|\K_{s,t}\| &\le& C\int_s^t \|\J_{r,t}\| \|\J_{s,r}\|_1^{3/2}
\|\J_{s,r}\|^{1/2}\,dr\\ &\le& C \exp \Bigl(\eta \int_0^1
\|w_r\|_1^2\,dr \Bigr) \;,
\end{eqnarray*}
where we used the integrability of $|r-s|^{-3/4}$ in the second step.
This concludes the proof of \lem{lem:mainbounds}.
\hfq

\demo{Proof of Lemma {\rm \ref{l:HighContract}}} 
  In order to get \eref{e:boundPhJ1}, we show that with the above
  notation, one can get bounds on $\|\uu_t\|_1$ as well. To achieve this
  we define, for a constant $\eps > 0$ to be fixed later, $\zeta_t =
  \|\uu_t\|^2 + t \eps \|\uu_t\|_1^2$. Using \eref{e:boundxi1} to bound
  the derivative of the first term and combining \eref{e:boundB} with
  \eref{e:idK} for the other terms, we then get in a straightforward
  way
\begin{eqnarray*}
  \d_t \zeta_t &\le& (\eps-\nu) \|\nabla \uu_t\|^2 + {C \over \eta^2
    \nu} \|\uu_t\|^2 + {\eta\over 2} \|w_t\|_1^2 \|\uu_t\|^2\\ &&
  -2t \eps \nu \|\uu_t\|_2^2 + 2t\eps C\|w_t\|_1
  \|\uu_t\|_2\|\uu_t\|_{1/2}\;.
\end{eqnarray*}
By \eref{e:interp}, we get
\begin{eqnarray*}
  2 C \|w\|_1 \|\uu\|_2\|\uu\|_{1/2} &\le& 2\sqrt{\eta\nu} \|w\|_1
  \|\uu\|_2\|\uu\|_1 + {C\over \eta\nu} \|w\|_1 \|\uu\|_2 \|\uu\| \\
  &\le& \nu \|\uu\|_2^2 + \eta \|w\|_1^2 \|\uu\|_1^2 + {C \over \eta^2
    \nu^3}\|w\|_1^2 \|\uu\|^2\;.
\end{eqnarray*}
This immediately yields
\begin{eqnarray*}
  \d_t \zeta_t &\le& (\eps-\nu) \|\uu_t\|_1^2 + {C \over \eta^2 \nu}
  \|\uu_t\|^2 + \Bigl({\eta\over 2} + {t\eps C \over \eta^2
    \nu^3}\Bigr) \|w_t\|_1^2 \|\uu_t\|^2\\ && - t \eps \nu
  \|\uu_t\|_2^2 + t\eps \eta \|w_t\|_1^2 \|\uu_t\|_1^2\;.
\end{eqnarray*}
If we take $\eps$ sufficiently small (of the order $\eta^3\nu^3$), we
get
$$
  \d_t \zeta_t \le \Bigl({C \over \eta^2 \nu} + \eta \|w_t\|_1^2\Bigr)
  \zeta_t\;,\quad \text{for $t \in [0,1]$,}
$$   
and therefore
\begin{equation} \label{e:boundzeta}
\|\J_t \uu_0\|_1^2 \le {\tilde C\over t} \exp \Bigl(\eta \int_0^1
\|w_s\|_1^2\,ds \Bigr) \|\uu_0\|^2\;,
\end{equation}   
for some (possibly rather large) constant $\tilde C$.  If we now
define $\pi_N$ as the orthogonal projection on the set of Fourier modes
with $|k| \ge N$,  
$$
\|\pi_N\uu_t\| \le {1\over N}\|\uu_t\|_1\;.
$$   
The bound \eref{e:boundPhJ1} immediately follows by taking $\pi_\ell = 1-\pi_N$
for $N$ sufficiently large.

We now turn to the proof of the bound \eref{e:boundPhJ2}. We define
$\pi_\ell$ as above (but reserve the right to choose the precise value
of $N$ later) and set $\uu_t^\ell = \pi_\ell \uu_t$ and $\uu_t^h =
(1-\pi_\ell)\uu_t$. With this notation, \eref{e:boundPhJ2} amounts
to obtaining bounds on $\|\uu_t\|$ with $\uu_0^\ell = 0$. Using the identity
\eref{e:idB} we have
\begin{eqnarray*}
  \d_t \|\uu_t^\ell\|^2 &=& -2\nu \| \uu_t^\ell\|_1^2 + 2
  \scal{B(\KK \uu_t^\ell, w_t), \uu_t^\ell} \\
&&- 2\scal{B(\KK \uu_t^h, \uu_t^\ell), w_t} 
   - 2 \scal{B(\KK w_t, \uu_t^\ell), \uu_t}\;, \\
  \d_t \|\uu_t^h\|^2 &= &-2\nu \|\uu_t^h\|_1^2 - 2 \scal{B(\KK
    \uu_t, \uu_t^h), w_t} - 2 \scal{B(\KK w_t, \uu_t^h), \uu_t}\;.
\end{eqnarray*}
Applying \eref{e:boundB} to the right-hand side allows us  to get the bound
\begin{eqnarray*}
  \d_t \|\uu_t^\ell\|^2 &\le &-2\nu \|\uu_t^\ell\|_1^2 + 2
  \scal{B(\KK \uu_t^\ell, w_t), \uu_t^\ell} + C\|w_t\|_{1/2} \|\uu_t^\ell\|_1 \|\uu_t^h\| \;, \\
  \d_t \|\uu_t^h\|^2 &\le& -2\nu \|\uu_t^h\|_1^2 + C\|w_t\|_{1/2} \|\uu_t^h\|_1 \|\uu_t\|\;.
\end{eqnarray*}
We then bound the first line using \eref{e:boundinterp}
and the second line using $\|\uu_t^h\|_1^2
\ge N^2 \|\uu_t^h\|^2$. We thus obtain
\begin{eqnarray}
  \d_t \|\uu_t^\ell\|^2 &\le &\Bigl({C\over \eta^2} + \eta \|w_t\|_{1}^2 
  \Bigr)\|\uu_t^\ell\|^2
  + C \|w_t\|_{1/2}^2 \|\uu_t^h\|^2\;, \label{e:boundlow}\\[4pt]
  \d_t \|\uu_t^h\|^2 &\le& -\nu N^2 \|\uu_t^h\|^2 + C \|w_t\|_{1/2}^2
  \|\uu_t\|^2\;, \nonumber
\end{eqnarray}
for an arbitrary value of $\eta$ and for a constant $C$ depending on
$\nu$ but independent of $N$ and $\eta$.  Using the \textit{a priori} bound
from point~3 above for the Jacobian $\uu$ and the interpolation inequality $\|w_s\|_{1/2}^2 \leq \|w_s\|\|w_s\|_{1}$
 immediately produces the bound
\begin{eqnarray*}
  \|\uu_t^h\|^2 &\le &e^{-\nu N^2 t} \|\uu_0^h\|^2 +  C \int_0^t
  e^{-\nu N^2 (t-s)}  \|w_s\|_{1/2}^2  \|\uu_s\|^2\,ds \\ [4pt]
  &\le &e^{-\nu N^2 t} \|\uu_0^h\|^2 + {C(T) \|\uu_0^h\|^2 \over N} e^{\eta
  \int_0^t \|w_s\|_1^2\,ds}\sqrt{\int_0^t \|w_s\|_1^2\,ds} \sup_{s\in[0,t]}\|w_s\|
  \\[4pt]
  &\le &\|\uu_0^h\|^2 \Bigl(e^{-\nu N^2 t} + {C(T) \over N}\exp \Bigl(\eta'
  \int_0^t \|w_s\|_1^2\,ds\Bigr)\sup_{s\in[0,t]}\|w_s\| \Bigr)\;,
\end{eqnarray*}
for an arbitrary $\eta' > \eta$.  Combining this
with the bound of point 1 above shows that, for every $\eta$, every
$\gamma$, every $p$, and every $T$, there exists a constant $N_0$ such
that
\begin{equation*}
  \expect_w \|\uu_t^h\|^p \le 2 e^{-\nu N^2 p t} \|\uu_0^h\|^p +
  \gamma e^{\eta \|w\|^2} \|\uu_0^h\|^p\;,
\end{equation*}
for all $t \in [0,T]$ and all $N \ge N_0$. Since $\uu_0^\ell = 0$ by
assumption, it follows from \eref{e:boundlow} that
\begin{eqnarray*}
  \|\uu_t^\ell\|^2 &\le& C \int_0^t \exp \Bigl({C(t-s)\over \eta^2}  +
  \eta \int_s^t \|w_r\|_1^2\,dr \Bigr) \|w_s\|_{1/2}^2 \|\uu_s^h\|^2\,ds\\ [4pt]
  &\le& \Bigl( \int_0^t e^{{8C(t-s)\over \eta^2} + 8\eta \int_s^t
    \|w_r\|_1^2\,dr}\,ds \Bigr)^{1/8} \Bigl(\int_0^t
  \|\uu_s^h\|^8\,ds\Bigr)^{1/4} \\[4pt]
  &&\times \sqrt{\int_0^t \|w_s\|_1^2\,ds} \sup_{s\in[0,t]}\|w_s\|\;.
\end{eqnarray*}
The required bound \eref{e:boundPhJ2} now follows easily by taking
expectations and using the previous bounds.
\Endproof\vskip12pt

\references {FJMR022}
 
  \def\url#1{\texttt{#1}} 

\bibitem[AN87]{MR88f:90180}
\name{E.~J. Anderson} and \name{P.~Nash},
 \emph{Linear programming in Infinite-Dimensional Spaces},
 {\it Wiley-Interscience Series in Discrete Mathematics and Optimization},
  John Wiley \& Sons Ltd., Chichester, U.K., 1987.

\bibitem[Bel87]{b:Bell87}
\name{D.~R. Bell},
 \emph{The {M}alliavin Calculus},
 Longman Scientific \& Technical, Harlow, 1987.

\bibitem[BKL01]{b:BricmontKupiainenLefevere01}
\name{J.~Bricmont}, \name{A.~Kupiainen}, and \name{R.~Lefevere},
 Ergodicity of the 2{D} {N}avier-{S}tokes equations with random
  forcing,
 \emph{Comm.\ Math.\ Phys.\/} \textbf{224}  (2001), 65--81.

\bibitem[BKL02]{b:BricmontKupiainenLefevere02}
\bibline,
 Exponential mixing of the 2{D} stochastic {N}avier-{S}tokes dynamics,
 \emph{Comm.\ Math. \ Phys.\/} \textbf{230}  (2002), 87--132.

\bibitem[Cer99]{CerrRD}
\name{S.~Cerrai},
 Ergodicity for stochastic reaction-diffusion systems with polynomial
  coefficients,
 \emph{Stochastics Stochastics Rep.\/} \textbf{67}  (1999),
  17--51.

\bibitem[CF88]{b:CoFo88}
\name{P.~Constantin} and \name{C.~Foia{\c{s}}},
 \emph{{N}avier-{S}tokes Equations},
 University of Chicago Press, Chicago, 1988.

\bibitem[CK97]{b:ChowKhasminskii98}
\name{P.-L. Chow} and \name{R.~Z. Khasminskii},
 Stationary solutions of nonlinear stochastic evolution equations,
 \emph{Stochastic Anal.\ Appl.\/} \textbf{15}  (1997), 671--699.

\bibitem[Cru89]{MR90c:35161}
\name{A.~B. Cruzeiro},
 Solutions et mesures invariantes pour des \'equations d'\'evolution
  stochastiques du type {N}avier-{S}tokes,
 \emph{Exposition. Math.\/} \textbf{7}  (1989), 73--82.

\bibitem[DPZ96]{b:DaZa96}
\name{G.~Da~Prato} and \name{J.~Zabczyk},
 \emph{Ergodicity for Infinite Dimensional Systems},
 Cambridge Univ.\ Press, Cambridge, 1996.

\bibitem[DM83]{DM83}  \name{C. Dellacherie}  and \name{P.-A. Meyer},
 \textit{Probabilit{\hskip.5pt\rm \'{\hskip-5pt\it e}}s et Potentiel. Chapitres {\rm IX}
{\rm \`{\hskip-5.5pt\it a}} {\rm XI:} Th{\hskip.5pt\rm \'{\hskip-5pt\it e}}orie discr{\hskip1pt\rm
\`{\hskip-5.5pt\it e}}te du Potentiel},
 {\it Actualit{\hskip.5pt\rm \'{\hskip-5pt\it e}}s Scientifiques et Industrielles} {\bf 1410},  Hermann, Paris
1983.

\bibitem[EH01]{EH3}
\name{J.-P. Eckmann} and \name{M.~Hairer},
 Uniqueness of the invariant measure for a stochastic {P}{D}{E} driven
  by degenerate noise,
 \emph{Comm.\  Math.\ Phys.\/} \textbf{219}  (2001),
  523--565.

\bibitem[EL94]{Xuemei}
\name{K.~D. Elworthy} and \name{X.-M. Li},
 Formulae for the derivatives of heat semigroups,\break
 \emph{J. Funct.\ Anal.\/} \textbf{125}  (1994), 252--286.

\bibitem[EM01]{b:EMattingly00}
\name{W.~E} and \name{J.~C. Mattingly},
 Ergodicity for the {N}avier-{S}tokes equation with degenerate random
  forcing: finite-dimensional approximation,
 \emph{Comm.\ Pure Appl.\ Math.\/} \textbf{54}  (2001),
  1386--1402.

\bibitem[EMS01]{b:EMattinglySinai00}
\name{W.~E. Mattingly}, \name{J.~C. Mattingly}, and \name{Y.~G. Sinai},
 Gibbsian dynamics and ergodicity for the stochastic forced
  {N}avier-{S}tokes equation,
 \emph{Comm.\ Math. \ Phys.\/} \textbf{224} (2001), 83--106.

\bibitem[Fer97]{b:Fe97}
\name{B.~Ferrario},
 Ergodic results for stochastic {N}avier-{S}tokes equation,
 \emph{Stochastics and Stochastics Reports} \textbf{60} 
  (1997), 271--288.

\bibitem[FJMR02]{MR1914189}
\name{C.~Foias}, \name{M.~S. Jolly}, \name{O.~P. Manley}, and
  \name{R.~Rosa},
 Statistical estimates for the {N}avier-{S}tokes equations and the
  {K}raichnan theory of 2-{D} fully developed turbulence,
 \emph{J. Statist.\ Phys.\/} \textbf{108}  (2002), 591--645.

\bibitem[Fla94]{b:Fl94}
\name{F.~Flandoli},
 Dissipativity and invariant measures for stochastic {N}avier-{S}tokes
  equations,
 \emph{Nonlinear Differential Equations Appl.\/} \textbf{1}  (1994), 403--426.

\bibitem[FM95]{FM}
\name{F.~Flandoli} and \name{B.~Maslowski},
 Ergodicity of the 2-{D} {N}avier-{S}tokes equation under random
  perturbations,
 \emph{Comm.\ Math.\ Phys.\/} \textbf{172} (1995),
  119--141.

\bibitem[FP67]{b:FoiasProdi67}
\name{C.~Foia{\c{s}}} and \name{G.~Prodi},
 Sur le comportement global des solutions nonstationnaires des
  \'equations de {N}avier-{S}tokes en dimension {$2$},
 \emph{Rend.\ Sem.\ Mat.\ Univ.\ Padova} \textbf{39} (1967), 1--34.

\bibitem[Fri95]{b:Frisch95}
\name{U.~Frisch},
 \emph{Turbulence\/{\rm :}\/ The Legacy of A. N. Kolmogorov}, Cambridge Univ.\ Press,
 Cambridge, 1995.

\bibitem[Hai02]{HExp02}
\name{M.~Hairer},
 Exponential mixing properties of stochastic {PDE}s through asymptotic
  coupling,
 \emph{Probab.\ Theory Related Fields} \textbf{124} (2002),
  345--380.

\bibitem[H{\"o}r67]{H1}
\name{L.~H{\"o}rmander},
 Hypoelliptic second order differential equations,
 \emph{Acta Math.\/} \textbf{119}  (1967), 147--171.

\bibitem[H{\"o}r85]{Ho}
\bibline
 \emph{The Analysis of Linear Partial Differential Operators {\rm 
  1--\rm  4}},
 {Springer-Verlag}, {New York}, 1985.

\bibitem[Kan42]{Kantorovich}
\name{L.~V. Kantorovich},
 On the translocation of masses,
 \emph{Dokl.\ Akad.\ Nauk SSSR} \textbf{37}  (1942), 194--201.

\bibitem[Kan48]{Kantorovich2}
\bibline,
 On a problem of {M}onge,
 \emph{Uspekhi Mat.\ Nauk} \textbf{3}  (1948), 225--226.

\bibitem[KS00]{b:KuksinShirikyan00}
\name{S.~Kuksin} and \name{A.~Shirikyan},
 Stochastic dissipative {P}{D}{E}s and {G}ibbs measures,
 \emph{Comm.\ Math. \ Phys.\/} \textbf{213}  (2000), 291--330.

\bibitem[KS01]{ArmenKuk}
\bibline,
 A coupling approach to randomly forced nonlinear {PDE}'s.~{I},
 \emph{Comm.\  Math.\  Phys.\/} \textbf{221}  (2001), 351--366.

\bibitem[Mat98]{b:Mattingly98b}
\name{J.~C. Mattingly}.
 {The stochastically forced {N}avier-{S}tokes equations: energy
  estimates and phase space contraction},
 Ph.D. thesis, Princeton University, 1998.

\bibitem[Mat99]{b:Mattingly98}
\bibline,
 Ergodicity of $2${D} {N}avier-{S}tokes equations with random forcing
  and large viscosity,
 \emph{Comm.\ Math. \ Phys.\/} \textbf{206}  (1999), 273--288.

\bibitem[Mat02a]{b:Mattingly02d}
\bibline,
 The dissipative scale of the stochastics {N}avier-{S}tokes equation:
  regularization and analyticity,
 \emph{J. Statist.\ Phys.\/} \textbf{108}  (2002), 1157--1179.

\bibitem[Mat02b]{MatNS}
\bibline,
 Exponential convergence for the stochastically forced
  {N}avier-{S}tokes equations and other partially dissipative dynamics,
 \emph{Comm.\ Math.\ Phys.\/} \textbf{230}  (2002),
  421--462.

\bibitem[Mat03]{b:Mattingly03Pre}
\bibline,
 On recent progress for the stochastic {N}avier-{S}tokes equations,
 In\break \emph{Journ{\hskip.5pt\rm \'{\hskip-5pt\it e}}es \'Equations aux
 d{\hskip.5pt\rm \'{\hskip-5pt\it e}}riv{\hskip.5pt\rm \'{\hskip-5pt\it e}}es partielles},
  Forges-les-Eaux, 2003.

\bibitem[MP06]{b:MattinglyPardoux03Pre}
\name{J.~C. Mattingly} and \name{E.~Pardoux}.
 Malliavin calculus and the randomly forced {N}avier-{S}tokes
  equation, {\it Comm.\ Pure Appl.\ Math.\/} {\bf 59} (2006), 1742--1790; math.PR/0407215.

\bibitem[MS05]{b:MattinglySuidan05}
\name{J.~C. Mattingly} and \name{T.~M. Suidan},
 The small scales of the stochastic {N}avier-{S}tokes equations under
  rough forcing,
 \emph{J. Stat.\ Phys.\/} \textbf{118}  (2005), 343--364.

\bibitem[MY02]{b:MasmoudiYoung02}
\name{N.~Masmoudi} and \name{L.-S. Young},
 Ergodic theory of infinite dimensional systems with applications to
  dissipative parabolic {PDE}s,
 \emph{Comm.\ Math. \ Phys.\/} \textbf{227} (2002), 461--481.

\bibitem[MR04]{b:MikuleviciusRozovskii02}
\name{R.~Mikulevicius} and \name{B.~L. Rozovskii},
 Stochastic {N}avier-{S}tokes equations for turbulent flows,
 \emph{SIAM J. Math.\ Anal.\/} \textbf{35}  (2004), 1250--1310
  (electronic).

\bibitem[Nua95]{MR96k:60130}
\name{D.~Nualart},
 \emph{The {M}alliavin calculus and related topics},
 {\it Probability and its Applications}, Springer-Verlag, New
  York, 1995.

\bibitem[Rac91]{MR93b:60012}
\name{S.~T. Rachev},
 {Probability metrics and the stability of stochastic models},
 {\it Wiley Series in Probability and Mathematical Statistics\/{\rm :}\/ Applied
  Probability and Statistics}, John Wiley \& Sons Ltd., Chichester, U.K., 1991.

\bibitem[Ros02]{MR2003k:76071}
\name{R.~M.~S. Rosa},
 Some results on the {N}avier-{S}tokes equations in connection with
  the statistical theory of stationary turbulence, in  {\it Mathematical Theory in Fluid Mechanics}
(Paseky, 2001),  \emph{Appl.\ Math.\/} \textbf{47}  (2002), 485--516.
 
\bibitem[Sei02]{Sei02} \name{J. Seidler}, {A note on the strong Feller property},
unpublished lecture notes  (2001), \texttt{http://simu0292.utia.cas.cz/seidler/nsfp.ps}.

\Endrefs

\end{document}